\documentclass{article}

\usepackage{amssymb,amsfonts}
\usepackage[all,arc]{xy}
\usepackage{enumerate}
\usepackage{mathrsfs}
\usepackage[T1]{fontenc}
\usepackage[utf8]{inputenc}
\usepackage[english,activeacute]{babel}
\usepackage{amssymb,amsmath,amsfonts,eucal,xspace,enumerate,amsthm,mathtools}
\usepackage{bbm}
\usepackage{cite}
\usepackage{graphicx}
\usepackage{float}
\usepackage{siunitx}
\usepackage{parskip}
\usepackage[font=small,labelfont=bf]{caption}

\usepackage{tikz, caption}

\swapnumbers\newtheorem{thm}{Theorem}[section]
\newtheorem{defi}[thm]{Definition}
\newtheorem{conj}[thm]{Conjecture}
\newtheorem{notation}[thm]{Notation}

\newtheorem{rem}[thm]{Remark}
\newtheorem{exam}[thm]{Example}
\newtheorem{lem}[thm]{Lemma}
\newtheorem{prop}[thm]{Proposition}
\newtheorem{coro}[thm]{Corollary}
\newtheorem{obs}[thm]{Observation}

\newcommand{\bobs}{\begin{obs}}
\newcommand{\eobs}{\end{obs}}

\newcommand{\blem}{\begin{lem}}
\newcommand{\elem}{\end{lem}}

\newcommand{\bteo}{\begin{thm}}
\newcommand{\eteo}{\end{thm}}
\newcommand{\bpro}{\begin{prop}}
\newcommand{\epro}{\end{prop}}
\newcommand{\bcoro}{\begin{coro}}
\newcommand{\ecoro}{\end{coro}}
\newcommand{\brem}{\begin{rem}}
\newcommand{\erem}{\end{rem}}
\newcommand{\bexa}{\begin{exam}}
\newcommand{\eexa}{\end{exam}}
\newcommand{\bdn}{\begin{defi}}
\newcommand{\edn}{\end{defi}}

\newcommand{\bdem}{\begin{proof}}
\newcommand{\edem}{\end{proof}}

\newcommand{\bconj}{\begin{conj}}
\newcommand{\econj}{\end{conj}}

\newcommand{\tbdem}{\begin{proof}[\proofname\ Theorem~\ref{Poly decay Neutral geom Lorenz flow Thm}]}
\newcommand{\tedem}{\end{proof}}

\newcommand{\bnotes}{\begin{notation}}
\newcommand{\enotes}{\end{notation}}

\newcommand{\bdes}{\begin{description}}
\newcommand{\edes}{\end{description}}

\newcommand{\tr}{\operatorname{trace}}

\title{\textbf{Mixing rates of the geometrical neutral Lorenz model}}

\author{Henk Bruin$^1$ \and Hector Homero \\Canales Farías$^2$} 
\date{%
    $^1$\textit{Faculty of Mathematics, University of Vienna}\\%
    $^2$\textit{Faculty of Mathematics, University of Vienna}\\[2ex]
    \today}

\newcommand{\bigo}[1]{\mathcal{O}(#1)}

\newcommand{\R}{\mathbb{R}}

\newcommand{\Cdiff}[1]{\mathit{C}^{#1}}

\newcommand{\Z}{\mathbb{Z}}

\newcommand{\fun}[3]{#1 : #2 \to #3}

\newcommand{\Div}[1]{\operatorname{Div}(#1)}

\newcommand{\norm}[2]{\left\lVert #1\right\rVert_{#2}}

\newcommand{\lp}[2]{\operatorname{L}^{#1}(#2)}

\newcommand{\on}[1]{\operatorname{#1}}
\newcommand{\vertiii}[2]{{\left\vert\kern-0.25ex\left\vert\kern-0.25ex\left\vert #1 
    \right\vert\kern-0.25ex\right\vert\kern-0.25ex\right\vert}_{#2}}
 
\newcommand{\ex}[1]{\operatorname{e}^{#1}}
\newcommand{\vf}[2]{\mathfrak{#1}^{#2}}
\newcommand{\fneu}{f_{\on{Neu}}}
\newcommand{\fLor}{f_{\on{Lor}}}
\newcommand{\gneu}{g_{\on{Neu}}}
\newcommand{\gLor}{g_{\on{Lor}}}
\newcommand{\Pneu}{P_{\on{Neu}}}
\newcommand{\Plor}{P_{\on{Lor}}}

\newcommand{\tnxbt}[1]{\tau'_{\on{Neu}}(#1)=\Abtwo{#1}{\beta_2}}

\newcommand{\rlor}[1]{r_{\on{Lor}}(#1)}
\newcommand{\rneu}[1]{r_{\on{Neu}}(#1)}

\newcommand{\Ab}[2]{A_1(#1,#2)}
\newcommand{\Abtwo}[2]{A_2(#1,#2)}
\newcommand{\cW}{\mathcal W}

\begin{document}
\maketitle
\tableofcontents

\begin{abstract}
The aim of this paper is to obtain polynomial decay of correlations
of a Lorenz-like flow where the hyperbolic saddle at the origin is replaced by a neutral saddle.
To do that, we take the construction of the geometrical Lorenz flow and proceed by changing the nature of the saddle fixed point at the origin by a neutral fixed point. This modification is accomplished by changing the linearised vector field in a neighbourhood of the origin for a neutral vector field.
This change in the nature of the fixed point will produce polynomial tails for the Dulac times,
and combined with methods of Ara\'ujo and Melbourne (used to prove exponential mixing for the classical Lorenz flow)
this will ultimately lead to polynomial upper
bounds of the decay of correlations for the modified flow.
\\[3mm]
\\[2mm]
\end{abstract}

\section{Introduction}\label{Intro}

The study of flows on surfaces and higher-dimensional manifolds has caught the interest of many scientists because of its numerous applications such as Hamiltonian flows, geodesic and horocycle flows, billiard flows or flows from meteorological models. These flows are usually equipped with a natural invariant measure $\mu$, for instance the SRB-measure.

The main goal is to have a better understanding of the properties of these flows, such as hyperbolicity, ergodicity, mixing (or at least  weak mixing) and, in chaotic settings, rates of mixing; that is, we would like to investigate the asymptotic behaviour of the correlation coefficients

\begin{equation}
\rho_t(v,w)=\Bigg|\int_M v\cdot w\circ f^td\mu-\int_M vd\mu\int_M wd\mu\Bigg|,
\end{equation}
where $\fun{f^t}{M}{M}$ is a flow acting on a manifold $M$ and $\mu$ its SRB-measure, and for observables $v$, $w$ chosen from an appropriate Banach space.
Knowing the rates of mixing is very helpful for proving other ergodic and statistical properties since mixing is one of the strongest statistical properties.

Obtaining good mixing rates for flows, even for hyperbolic flows, is far more difficult than for maps. Some seminal ideas were provided by Liverani \cite{Liv04} and Dolgopyat \cite{DolAnosov,DolDio}, with applications of these methods in e.g.\ \cite{ExpMixTeich, BB05, BT18}.
To obtain sharp estimates in the polynomial setting, the operator renewal theory techniques developed by Sarig \cite{Sarig} and Gou\"ezel \cite{Gouezel} are the only ones available.

The model we would like to study is probably one of the most emblematic ones, the Lorenz flow. In the mid seventies
 Afra\u{\i}movi\v{c}, Bykov and Shilnikov
\cite{ABS77} and independently Guckenheimer and Williams \cite{GW}
introduced the geometric Lorenz attractor to model the original Lorenz attractor.
Our research focuses on a modified version of this geometrical model
and study its rate of mixing, based on the precise estimates
in \cite{Regular} of Dulac times associated to a neutral saddle.

Recently, Ara{\'u}jo and Melbourne in \cite{ExpDecayNonUniform} proved that the geometrical Lorenz flow 
(and hence the classical Lorenz flow), also enjoys exponential mixing. 
It is techniques from their papers, specifically $C^{1+\alpha}$ smoothness of the stable foliation,
that leads eventually to the claimed mixing rates.

\subsection{The framework}

The geometrical Lorenz flow can be seen as the natural extension of a suspension 
semiflow built over a certain type of one-dimensional expanding map $\fLor$. We first consider the cross-section $\Sigma$ transversal to the flow and the Poincar\'e map $\fun{\Plor}{\Sigma}{\Sigma}$, which is decomposed in two parts. The first one is the Dulac map, denoted by $P_1$, deals with the local behaviour near the origin and is obtained by considering a linear system in a neighbourhood of the origin; that is,
we take the flow $\on{X}^t$ obtained from the linear system

\begin{equation}\label{LorenzLinearized} 
\begin{array}{ccc}
\dot{x}=\\
\dot{y}=\\
\dot{z}=\\
\end{array}
\begin{array}{ccc}
\lambda_u x\\
-\lambda_s y\\
-\lambda_{ss} z\\
\end{array}
\end{equation}
where $\lambda_u$, $\lambda_s$ and $\lambda_{ss}$ denote the unstable, stable and strong stable
eigenvalues of the original Lorenz system, respectively. Then we let points in $\Sigma$ flow under $\on{X}^t$ 
until flow time $\tau'_{\on{Lor}} := \min\{ t > 0 : \on{X}^t \in S\} = -\lambda_{u}^{-1}\ln(|x|)+\mathcal{O}(\ln(|x|))$ as $x \to 0$. Thus we have that $\fun{X^{\tau'_{\on{Lor}}}=P_1}{\Sigma}{S}$,
where $S^{\pm}$ is the image of $\Sigma^{\pm}$ under $P_1$ and has a cusp-like shape, see Figure~\ref{fig:Lorenz}.

The second part, denoted by $P_2$, consists of the return of $S$ to $\Sigma$ and simulates the random turns of a 
regular orbit around the origin and describes a butterfly-like shape. 
This is done by a composition of a rotation, expansion and translation with hitting time $\tau_2(x)\in\Cdiff{\epsilon}$. Thus, the full return time of the Poincar\'e map $\Plor=P_2\circ P_1$ is given by

\begin{equation}
\label{EqReturnTimeLor}
\rlor{x}=\tau'_{\on{Lor}}(x)+\tau_2(x).
\end{equation}

\begin{center}\vspace{1cm}
\includegraphics[width=0.85\linewidth]{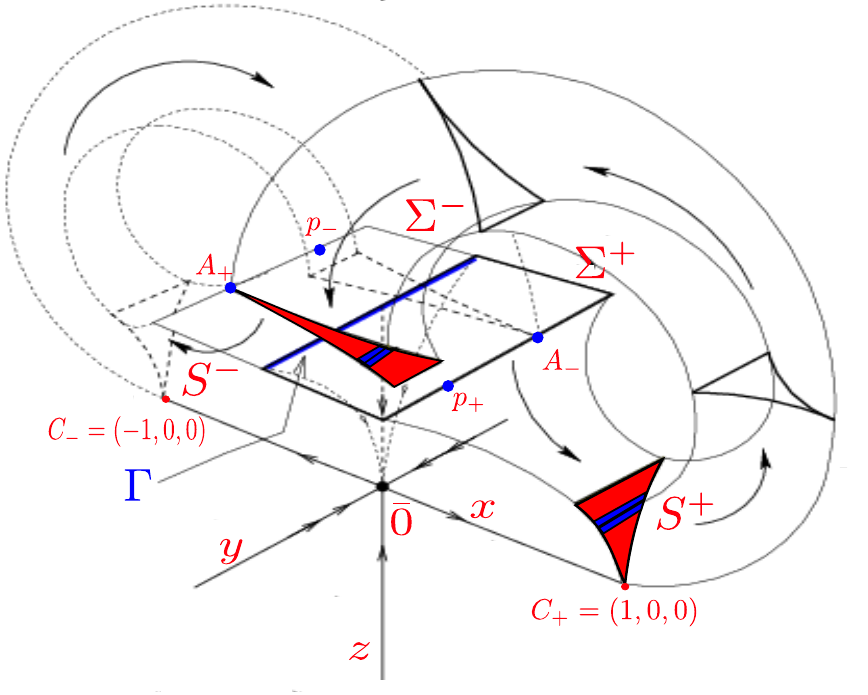}
\captionof{figure}{The Guckenheimer-Williams geometric model of the Lorenz flow (Image taken from \cite{3DimFlow})}
\label{fig:Lorenz}
\end{center}\vspace{1cm}

As we will see later, the lines in the $y$-direction (\textit{i.e.,},  parallel to the $y$ axis) in $\Sigma$ form the stable foliation, invariant under $P_{Lor}$; that is, for any leaf $\gamma$ of this foliation, its image $\Plor(\gamma)$ is contained in a leaf of the same foliation. By quotienting out the stable direction we can rewrite the Poincar\'e map as a skew-product; that is, $\Plor(x,y)=(\fLor(x),\gLor(x,y))$.

The \textbf{geometric Lorenz flow} is the couple $(W,X^t_W)$, where $W=\{\on{X}^t(\bar{x})\;|\;\bar{x}\in\Sigma,\;t\in\R^+\}$.
Consider $U=\bigcup_{\bar{x}\in\Sigma}X_W^{[0,\rlor{x}]}(\bar{x})$, then the \textbf{geometric Lorenz attractor} (of the corresponding vector field) is given by $\Lambda_{\on{Lor}}=\bigcap_{t>0} X_W^t(U)$.

In \cite{ExpDecayNonUniform}, exponential mixing for the geometrical Lorenz flow was proven under two conditions: the stable foliation has to be $\Cdiff{1+\alpha}$ and a uniform non-integrability (UNI) 
condition needs to be satisfied.

The modified version is obtained by changing the local behaviour near the origin. We achieve this by replacing the linear system for the following system;

\begin{equation} \label{Neutral system} 
\left( \begin{array}{ccc}
\dot{x}\\
\dot{y}\\
\dot{z}\\
\end{array}\right)=
Z\left( \begin{array}{ccc}
x\\
y\\
z\\
\end{array}\right)=
\left( \begin{array}{ccc}
x(a_0x^2 + a_1y^2 + a_2z^2)\\
-\ell y(1 + c_0x^2 + c_2z^2)\\
-z(b_0x^2 + b_1y^2 + b_2z^2)\\
\end{array}\right)+\mathcal{O}(4), 
\end{equation}
where $a_0,\;a_1,\;a_2,\;b_0,\;b_1,\;b_2,\;c_0,\;c_2$ and $\ell>0$,
$a_2b_0 < 9a_0b_2$, $\Delta:=a_2b_0-a_0b_2\neq0$
and $\mathcal{O}(4)$ refers to terms of order four or higher, under the condition that they are of the 
form $x^2\mathcal{O}(2)$ near the $yz$-plane and $z^2\mathcal{O}(2)$ near the $xy$-plane. 
This system has a polynomial Dulac time (see~\eqref{Dulac Map formula} and Figure~\ref{fig:Dulac map}) given by;

\begin{equation}\label{pol flow time N1}
\tau'_{\on{Neu}} := \min\{ t > 0 : \on{N}^t \in S\} = |x|^{-\frac{1}{\beta_2}}(1+\mathcal{O}(|x|^{\frac{1}{2\beta_2}})),
\end{equation}

\noindent as $x\to0$ and $\beta_2=\frac{a_2+b_2}{2b_2}$. To obtain the flow time $\tau'_{Neu}$, we make use of the estimates of the Dulac map and the tails of the return map obtained by Bruin and Terhesiu in \cite{Regular}. This change of flow time, from logarithmic to polynomial, 
will ultimately allow us to deduce the polynomial decay of correlations.

We denote by $\on{N}^t$ the flow obtained from the system given by~\eqref{Neutral system}. This change in the local behaviour near the origin leads to a change on the map $P_1$; that is, we have now $\fun{\on{N}^{\tau'_{\on{Neu}}}=D_1}{\Sigma}{S}$. For the second part, the return of $S$ to $\Sigma$, we consider the same diffeomorphism $P_2$ with same hitting time. In this way, we obtained the modified Poincar\'e map $\Pneu=P_2\circ d_1$ with return time given by,

\begin{equation}
\label{EqReturnTimeNeu}
r_{\on{Neu}}(x)=\tau'_{\on{Neu}}(x)+\tau_2(x).
\end{equation} 

Similarly, we define the \textbf{geometric neutral Lorenz flow} as the couple \\$(W,\on{N}^t_W)$, where $W=\{\on{N}^t(\bar{x})\;|\;\bar{x}\in\Sigma,\;t\in\R^+\}$. We consider again $U=\bigcup_{\bar{x}\in\Sigma}\on{N}^{[0,\rneu{x}]}(\bar{x})$, the \textbf{geometric neutral Lorenz attractor} (of the corresponding vector field) is given by $\Lambda_{\on{Neu}}=\bigcap_{t>0}\on{N}^t(U)$.

As we will see in more detail in Section \ref{Sec. Poincare maps Neutral models}, the geometrical neutral Lorenz flow will be split into three models. {\bf Model 1} is obtained when we take the parameters $c_0=c_2=0$  in \eqref{Neutral system}. {\bf Model 2} when we consider  $a_1=b_1=0$. Finally, {\bf Model 3}, the most general, will be given by taking all parameters strictly positive.

\subsection{Main results}

Let $C^{\eta}$ be the space of  functions that are $\eta$-H\"older
in the space direction, and  $C^{m,\eta}$ be the space of functions that are $m+\eta$-H\"older (i.e., $m$ time differentiable with an $\eta$-H\"older $m$-th derivative)
in the flow direction, see Section~\ref{Sec. proof of decay of correlations Neutral}
for the precise definitions.
The main result in this paper is the following theorem:


\bteo\label{Poly decay Neutral geom Lorenz flow Thm}
Let  $\fun{\on{N}^t}{\Lambda_{\on{Neu}}}{\Lambda_{\on{Neu}}}$ be the geometrical neutral Lorenz flow for Model 1 and Model 2 obtained from the neutral form given by \eqref{Neutral system}, with corresponding parameters. $\Lambda_{\on{Neu}}$ its attractor and its SRB measure $\mu$.
Then $N^t$ has polynomial decay of correlations (with exponent $\beta_2=\frac{a_2+b_2}{2b_2}$); that is, there exist $m\geq1$ and a constant $C>0$ such that for observables $v\in C^{\eta}(M)\cap C^{0,\eta}(M)$, $w\in C^{m,\eta}(M)$, and time $t>1$ we have
\begin{align*}
\rho_t(v,w)\leq & C(\norm{v}{{C^{\eta}}}+\norm{v}{{C^{0,\eta}}})\norm{w}{{C^{m,\eta}}}  t^{-\beta_2}.
\end{align*}
\eteo

A first question that presents itself is of course if these bounds are sharp,
and if current operator renewal theory methods \cite{Sarig, Gouezel}  cannot prove that.
We say more on this at the end of Section~\ref{Sec. proof of decay of correlations Neutral}.


For the proof of Theorem~\ref{Poly decay Neutral geom Lorenz flow Thm}, we obtain an explicit form of the Poincar\'e map, since we can solve the differential equation in the $y$ component. Thus we are able to prove polynomial decay of correlations by using the results on non-uniformly hyperbolic flows established by B\'alint \textit{et al.} in \cite{PolynomialDecayBalint}. 

For the third model the situation is more subtle since, to our knowledge, finding the solution of the differential equation in the $y$ component is next to impossible. To overcome this problem we will analyse and compare, with numerical methods, the limit behaviour of the Dulac maps obtained in \cite{Volume} and \cite{Regular} and adapted to our framework. More precisely, we will analyse the limit behaviour of the maps $\fun{D_1}{\Sigma}{S}$ obtained for each Neutral model. This is sufficient since the Poincar\'e maps considered in this work are given by $\Pneu=P_2\circ D_1$,
where $P_2$ is a diffeomorphism and the map $D_1$ is the Dulac map from the cross-section $\Sigma$ to the cusps $S$, which depends on the differential equation being considered. Therefore, the behavioural changes exhibited by the map $\Pneu$ are represented by the changes of the map $D_1$.
Dulac in \cite{Dulac} made a significant contribution to solving Hilbert's 16th problem by incorporating his map as an element to establish that polynomial vector fields in the plane possess a limited number of limit cycles, demonstrating that they cannot have an infinite number of such cycles.

The numerical analysis on the behaviour of the Dulac maps will give us  the plausibility of the following conjecture.

\bconj
Let  $\fun{\on{N}^t}{\Lambda_{\on{Neu}}}{\Lambda_{\on{Neu}}}$ be the geometrical neutral Lorenz flow for Model 3 obtained from the neutral form given by \eqref{Neutral system}, with the corresponding parameters. $\Lambda_{\on{Neu}}$ its attractor and its SRB measure $\mu$. Then $N^t$ has polynomial decay of correlations (with exponent $\beta_2=\frac{a_2+b_2}{2b_2}$).
\econj

The organization of this paper is as follows: In Section~\ref{Sec. Poincare maps Neutral models} 
we will give the construction of the Poincar\'e maps of the Neutral Model 1 and 2.
In Section~\ref{Sec. stable foliation and UNI} we will be devoted to the proof that the stable foliation for the geometrical neutral models is $\Cdiff{1+\alpha}$ and the UNI condition is satisfied by adapting the existing proofs for the geometrical Lorenz model. Section~\ref{Sec. proof of decay of correlations Neutral} contains the framework of non-uniformly hyperbolic flows and the proof of Theorem~\ref{Poly decay Neutral geom Lorenz flow Thm}. Finally, in Section~\ref{Sec. numerics} we will present the
numerical analysis and results we obtained for the Dulac map and the tails of the return map.    

For the remaining of this paper we will adopt the following notation.

\bnotes\label{notation higher terms}
In order to avoid excessive notation of the higher order terms, obtained from the estimates of the Dulac time given in \cite{Regular}, we will write $\Ab{x}{\beta}$ and $\Abtwo{x}{\beta_2}$ to denote  $\xi|x|^\beta(1+\bigo{|x|^{\frac{1}{2\beta_2}}})$ and $\zeta|x|^{-\frac{1}{\beta_2}}(1+\bigo{|x|^{\frac{1}{2\beta_2}}})$, respectively, where $\beta_0=\frac{a_0+b_0}{2a_0}$, $\beta_2=\frac{a_2+b_2}{2b_2}$, $\beta=\frac{\beta_0}{\beta_2}$, $\xi$ and $\zeta$ are constants given in \cite{Regular}, namely in Theorem 1.1 and the proof of Proposition 2.1.
$X\in\vf{X}{r}(M)$ will denote the vector space of $C^r$ vector fields in a manifold $M$ with the $C^r$ topology.
\enotes

{\em Acknowledgements:} This paper grew out of the thesis work of the second author. We gratefully acknowledge the support of FWF stand-alone grant P31950-N35.
We thank the referees for their very useful comments.

{\em Data availability statement:} The methods of how the numerical graphics were computed are given in the last section of this paper. For further details, please contact the authors.

\section{The Poincar\'e maps}\label{Sec. Poincare maps Neutral models}

\subsection{Neutral model 1} \label{Sec. model 1}

To create the modified models we will apply local surgery in a neighbourhood of the hyperbolic saddle equilibrium of 
the geometrical Lorenz model, namely the origin, and transform it into a neutral equilibrium. 
We do this because we aim to slow down the orbit and thus increase the time that orbits take to flow from the cross-section
$\Sigma$ to the cusps $S$, see Figure~\ref{fig:Lorenz},  and see the changes this new motion produces in the decay of correlations.
The flow obtained from this modification will be an almost Anosov flow \cite{ConditionsHu, Nonexistence}. Existence of a finite or infinite 
SRB measure for two-dimensional almost Anosov diffeomorphisms was already proven 
in \cite{ConditionsHu, Nonexistence}. Bruin and Terhesiu in \cite{Regular} proved mixing rates 
in the infinite SRB measure setting for almost Anosov diffeomorphism and established the
required spectral properties for the transfer operator (acting on an appropriate anisotropic Banach space
of distributions) of an induced map so as to obtain optimal rates of mixing. Furthermore, they gave more precise
tail estimates for the inducing scheme. We will take advantage of these methods and estimates 
and use them to deduce the rates of mixing of our almost Anosov flow.

We consider now $\Sigma^{*}=\{(x,y,1)\in\R^3\;|\;|x|\leq1,|y|\leq1\}$, $\Sigma^{-}=\{(x,y,1)\in\Sigma^{*}\;|\;x<0\}$,  $\Sigma^{+}=\{(x,y,1)\in\Sigma^{*}\;|\;x>0\}$, $\Sigma=\Sigma^{+}\cup\Sigma^{-}=\Sigma^{*}\setminus\tilde{\Gamma}$, where $\tilde{\Gamma}=\{(x,y,1)\in\Sigma^{*}\;|\;x=0,\}$
and $S=S^{+}\cup S^{-}$
ent where $S^{\pm}$ is the image of $\Sigma^{\pm}$ under the Dulac map $D_1$, see Figure~\ref{fig:Map N1}. The section $\Sigma$ is transversal to the flow and every trajectory eventually crosses $\Sigma$ in the direction of the negative axis $z$. Then for each $(x,y,1)\in\Sigma$, the time $\tau'_{\on{Neu}}$ such that $\on{N}^{\tau'_{\on{Neu}}}(x,y,1)\in S$ is determined by the estimates of the Dulac map provided in \cite{Regular}, as we will explain now.

\begin{center}\vspace{1cm}
\includegraphics[width=0.85\linewidth]{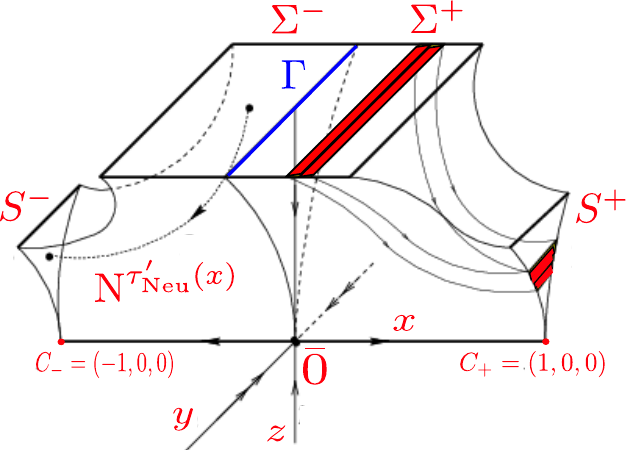}
\captionof{figure}{The map $D_1$ (Image taken from \cite{3DimFlow}).}
\label{fig:Map N1}
\end{center}\vspace{1cm}

Let us start with the Neutral model 1; that is, we consider a neighbourhood $U$ of the fixed point $\bar 0 = (0,0,0)$ where the vector field has the form \eqref{Neutral system} (in local Euclidean coordinates) with $c_0=c_2=0$ and the other parameters satisfying the constraints given before.
Note that $U$ is taken much smaller than the scale of Figure~\ref{fig:Map N1}.
This vector field, denoted by $Z$, is cubic at $\bar 0$ in the direction
transversal to the stable manifold of $\bar 0$, but this is the only source of non-hyperbolicity.
The $y$-axis is invariant and all solutions tend to $\bar 0$. The divergence is given by
$$
\Div{Z}=(3a_0-b_0)x^2+(a_1-b_1)y^2+(a_2-3b_2)z^2-\ell.
$$
Since we want the flow to shrink volume exponentially fast, as does the Lorenz flow,
we need $\Div{Z} \leq -c<0$.
Therefore, we let $\ell$ be large enough such that $(3a_0-b_0)x^2+(a_1-b_1)y^2+(a_2-3b_2)z^2<\ell$
for all $(x,y,z) \in U$. The solution for the $y$-component is given by $y(t)=y_0\ex{-\ell t}$. Thus we obtain a non-autonomous system of differential equations, since the contribution of the $y$-component to the $x$ and $z$- component is exponentially small as time increases, these terms are of smaller order than the higher order terms. Thus we are
left with the two-dimensional system studied in \cite{Regular}:

\begin{equation}
\label{FlowHor} 
\left( \begin{array}{ccc}
\dot{x}\\
\dot{z}\\
\end{array}\right)=
Z_{\on{hor}}\left( \begin{array}{ccc}
x\\
z\\
\end{array}\right)=
\left( \begin{array}{ccc}
x(a_0x^2+a_2z^2)\\
-z(b_0x^2+b_2z^2)\\
\end{array}\right).
\end{equation}

Now let $\cW^s$ and $\cW^u$ be two mutually tranversal foliations of
the interior of $U \cap \{\text{positive quadrant}\}$
that is invariant under the flow of \eqref{FlowHor}
and such that
\begin{itemize}
 \item the leaves of $\cW^s$ accumulate in $C^1$ topology on the stable manifold of $(0,0)$ and are transversal to the unstable manifold of $(0,0)$, and
 \item the leaves of $\cW^u$ accumulate in $C^1$ topology on the unstable manifold of $(0,0)$ and are transversal to the stable manifold of $(0,0)$.
\end{itemize}
One would like to use the stable and unstable foliation of the flow for $\cW^s$ and $\cW^u$, but as long as we defined the flow only locally, the above properties suffice.

Now fix an unstable leaf $W^u(0,z_0) \in \cW^s$ and a stable leaf $W^s(x_0,0) \in \cW^u$, then the Dulac map $\fun{D}{W^u(0,z_0)}{W^s(x_0,0)}$,
shown in Figure~\ref{fig:Dulac map}, assigns the first intersection $\phi^T(x,z_0)$ of the integral curve through $(x,z_0)$ with the stable leaf $W^s(x_0,0)$, where $x\in W^u(0,z_0)$, $\phi^t(x,z)$ is the flow from \eqref{FlowHor} and $T$ is the exit time. The estimates for the map $D$ and the flow time given in \cite{Regular} are:

\begin{equation}
\label{Dulac Map formula}
\omega=D(x)=c(z_0)x^\beta(1+\mathcal{O}(x^{\frac{1}{2\beta_2}}))
\end{equation}

\noindent and

\begin{equation}
\label{Map N1 flow time}
\tnxbt{x},
\end{equation}
where $\beta=\frac{\beta_0}{\beta_2}$ for $\beta_0=\frac{a_0+b_0}{2a_0}$ and $\beta_2=\frac{a_2+b_2}{2b_2}$.

\begin{center}\vspace{1cm}
\includegraphics[width=0.5\linewidth]{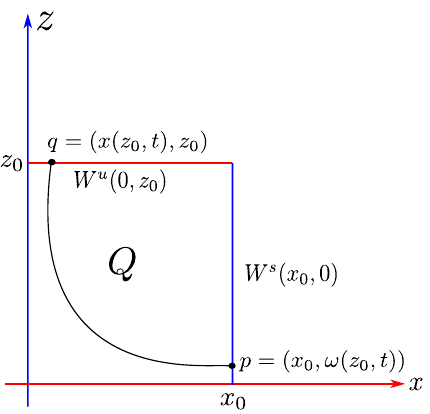}
\captionof{figure}{The Dulac map $\fun{D}{W^u(0,z_0)}{W^s(z_0,0)}$.}
\label{fig:Dulac map}
\end{center}\vspace{1cm}

More relevant to the proof of the decay of correlations of the neutral geometrical Lorenz flow is the estimate of the tails of the return map which we state in the following theorem.

\bteo \cite[Theorem 1.1]{Regular}\label{Estiamete tails of return} 
Let $\beta^{*}=\frac{1}{2}\min\{1,\frac{a_2}{b_2},\frac{b_0}{a_0}\}$, then there exists $C_0>0$ a constant such that
$\mu(\varphi>n)=C_0n^{-\beta_2}(1+\mathcal{O}(n^{-\beta^{*}}))$, where $\varphi=\min\{n\geq1\;|\; D_1^n(z)\notin P_0\}$ with $P_0$ an element of the Markov partition such that $Q = P_0 \cap \{\text{positive quadrant}\}$ (see Figure \ref{fig:Dulac map}) and $\mu$ the SRB-measure of $D_1^{\varphi}$.
\eteo

Putting all together we get the following expression for the map $D_1$ 

\begin{align}\label{Map N1 Neutral model1}
D_1(x,y,1) &=\on{N}^{\tau'_{\on{Neu}}(x)}(x,y,1) \nonumber \\
      &=\Big(1,y\ex{-\ell\Abtwo{x}{\beta_2}},\Ab{x}{\beta}\Big),
\end{align}
\noindent where the functions $\Abtwo{x}{\beta_2
}$ and $\Ab{x}{\beta}$ come from Notation \ref{notation higher terms}.

We make the following observations:

\bobs \label{Obs map N1} 
\begin{enumerate}
\item[1.-] $D_1(\Sigma^{\pm})$ has the shape of a cusp at $(\pm1,0,0)$ and (with some abuse of notation) we will denote these images as $S^{\pm}$ and $S=S^+\cup S^-$.

\item[2.-] Denote by $\ell_v(c)=\{(x,y,1)\in\Sigma\;|\;x=c\}$, where $c$ is a constant, the line segments in $\Sigma$ parallel to the $y$-axis and by $\ell_h(c)=\{(\pm1,y,z)\in S\;|\;z=c\}$, the line segments in $S$ parallel to the $y$-axis. Then  $D_1(\ell_v(c_0))=\ell_h(c_1)$; that is, the map $D_1$ takes the $y$-direction lines in $\Sigma$ to the horizontal line segments in $S$ as illustrated in Figure~\ref{fig:Map N1}.
\end{enumerate}
\eobs

The return of the cusps $S$ to the cross-section $\Sigma$ is described by the map $P_2=T\circ E_a\circ R_{\theta}$, where $R_\theta$ is a rotation by an angle of $\theta=\frac{3\pi}{2}$ and the rotation axis are the boundaries of the cross-section $\Sigma$ parallel to the $y$-axis, $E_a$ is an expansion by a factor of $a>1$ in the $x$-direction and a translation $T$ such that the unstable direction which starts from the origin is sent to the boundary of $\Sigma$; that is, we want to send the cusp points $C_{\pm}$ to $A_{\pm}$, see Figure~\ref{fig:Lorenz}. Thus, the full Poincar\'e map $\fun{\Pneu=P_2\circ D_1}{\Sigma}{\Sigma}$, shown in Figure~\ref{fig:Map PNeu}, is given by,

\begin{equation}
\label{P Neu Map1}
\Pneu(x,y) =
\begin{cases}
\left( a\Ab{x}{\beta}-1,\;y\ex{-\ell\Abtwo{x}{\beta_2}}-\dfrac{1}{2} \right),        & \text{if} \; x \in (0, 1];\\[1em]
\left( a\Ab{x}{\beta}+1,\;y\ex{-\ell\Abtwo{x}{\beta_2}}+\dfrac{1}{2}\right),        & \text{if} \; x \in [-1, 0).\\
\end{cases}
\end{equation}

In the positive quadrant the matrix $D\Pneu$ has eigenvalues $\lambda_1=a\frac{\beta_0}{\beta_2} |x|^{\frac{\beta_0}{\beta_2}-1}$ and $\lambda_2=\ex{-\ell\Abtwo{x}{\beta_2}}$. By restricting $\frac{1}{2}<\beta_0<2$ we have that $0<\dfrac{\beta_0}{\beta_2}<1$ since $\beta_2>2$. Then $\lambda_1\rightarrow\infty$ as $x\rightarrow0$. For the other eigenvalue we have that $\lambda_2<1$ and $\lambda_2\rightarrow0$ as $x \to 0$ since $\Abtwo{x}{\beta_2}\rightarrow\infty$ as $x \to 0$. Thus the modified Poincar\'e map $\Pneu$ is hyperbolic when $x$ approaches the origin.

\begin{center}\vspace{1cm}
\includegraphics[width=0.45\linewidth]{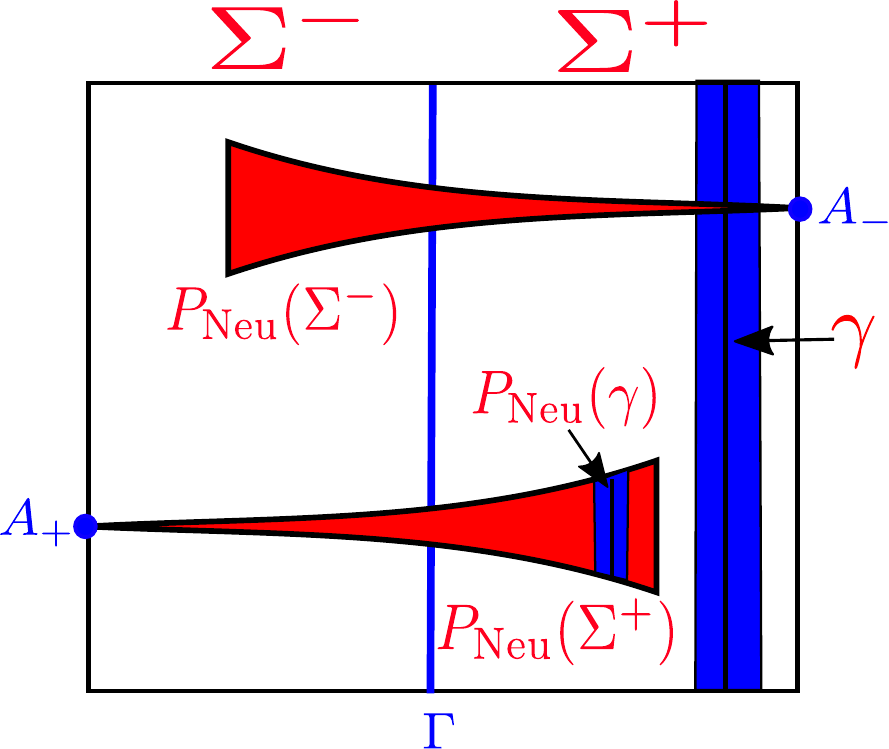}
\captionof{figure}{The Poincar\'e map for the neutral geometrical model.}
\label{fig:Map PNeu}
\end{center}\vspace{1cm}

The foliation given by the lines $\ell_v(c)$ is invariant under the map $\Pneu$; that is, given any leaf $\gamma$ of this foliation its image $\Pneu(\gamma)$ is contained in a leaf of the same foliation (See Figure~\ref{fig:Map PNeu}). Therefore, we can express $\Pneu$ as $\Pneu(x,y)=(\fneu(x),\gneu(x,y))$, where $\fun{\fneu}{I\setminus\{0\}}{I}$ is a Lorenz-like map with exponent $\beta=\frac{\beta_0}{\beta_2}$ with $\beta_0=\frac{a_0+b_0}{2a_0}$ and $\beta_2=\frac{a_2+b_2}{2b_2}$; that is, $\fneu$ is given by,

\begin{equation}\label{Lorenz map function}
\fneu(x)=
\begin{cases}
a\Ab{x}{\beta}-1,  & \text{if} \; x \in (0, 1];\\[1em]

1-a|\Abtwo{x}{\beta}|,  & \text{if} \; x \in [-1, 0),\\
\end{cases}
\end{equation}

\noindent and the function $\fun{\gneu}{(I\setminus\{0\})\times I}{I}$, where $I=[-1,1]$, satisfies the following:

\begin{enumerate}
\item[g1.-] The map $\gneu$ is piecewise $\Cdiff{2}$ and for fixed $x_0$, the map $\gneu(x_0,y)$ is a contraction in the $y$-direction, \textit{i.e.},
$$
d(\gneu(x_0,y_1),\gneu(x_0,y_2))\leq cd(y_1,y_2),
$$
where $d$ is the Euclidean distance in $I$ and $0<c<1$.

\item[g2.-] $D\Pneu$ has the following bound on its partial derivatives:

\begin{enumerate}
\item[a)] For all $(x,y)\in\Sigma$ we have $\partial_y\gneu(x,y)=\ex{-\ell\Abtwo{x}{\beta_2}}$. Since $\beta_2>2$ and $|x|\leq1$, there is $0<\eta<1$ such that

$$
|\partial_y\gneu(x,y)|<\eta.
$$

\item[b)] For $(x,y)\in\Sigma$ with $x\neq0$ we have $\partial_x\gneu(x,y)=\frac{\ell}{\beta_2}y(\Abtwo{x}{\beta_2})'\\\ex{-\ell\Abtwo{x}{\beta_2}}$. Since $1<1+\dfrac{1}{\beta_2}<\dfrac{3}{2}$ and $|y|,\;|x|\leq1$ we get that $|\partial_x\gneu(x,y)|$ is bounded. In fact, it tends to zero exponentially fast as $x$ approaches the origin.
\end{enumerate}

\item[g3.-] From g$2.$a) above follows the uniform contraction of the foliation given by the lines $\ell_v(c)$; in other words, there is $C>0$ such that, for any given leaf $\gamma$ of the foliation and for $y_1,\; y_2\in\gamma$, we have
$$
d(\Pneu^n(y_1),\Pneu^n(y_2))\leq C\eta^nd(y_1,y_2),
$$

\noindent when $n\rightarrow\infty$.
\end{enumerate}

\subsection{Neutral model 2} \label{Sec. model 2}

Now we will consider the Neutral model 2; that is, we consider the same neighbourhood $U$ of the origin where the flow has the local form given by \eqref{Neutral system} with $a_1=b_1=0$ and the remaining parameters satisfying the same constraints stated in this framework. Again $\bar 0=(0,0,0)$ is the only neutral periodic orbit and the vector field is cubic in the  direction transversal the stable manifold of $\bar 0$, but this is the only source of non-hyperbolicity. If $x=0$ and $z=0$, then we see that $\dot{x}=0$, $\dot{y}=-\ell y$ and $\dot{z}=0$. Hence, the $y$-axis is invariant and all solutions tend to the origin as in the previous model. Moreover, since $\dot{x}$ and $\dot{z}$ are decoupled from $y$, we have \eqref{FlowHor}. Thus the asymptotics for the Dulac map and the flow time given in \cite{Regular} follow. Also $\Div{Z}=(3a_0-\ell c_1-b_0)x^2+(a_2-\ell c_2-3b_2)z^2-\ell$.
Since we want a flow that shrinks volume exponentially fast as before, we take $\ell$ large enough so that $\frac{(3a_0-b_0)x^2+(a_2-3b_2)z^2}{1+c_0x^2+c_2z^2}<\ell$ for all $(x,y,z) \in U$.

We consider the same cross-section $\Sigma$ as before and proceed to construct the Poincar\'e map in the same way. Denote by $N^t$ the flow obtained from \eqref{Neutral system} with the pertinent constraints in the parameters; that is, $N^t(x,y,z)$ $=(x(t),y(t),z(t))$. By \eqref{Dulac Map formula} we obtain the following form for the flow,

\begin{equation}\label{NeutralFlowGeneralFormula}
N^{t(x,z)}(x,y,z)=(x_0,y(t(x,z)),\omega(z,t(x,z)).
\end{equation}
Note that $\dot{y}=y(-\ell(1+c_0x^2+c_2z^2))$, applying Gr\"onwall's Lemma we get,

\begin{align}
y(t)&=y_0\exp(-\ell\int_0^t(1+c_0x^2+c_2z^2)ds)\nonumber\\
&=y_0\ex{-\ell t}\exp(-\ell\int_0^t(c_0x^2+c_2z^2)ds).
\end{align}

By the estimates of the Dulac map and since $z\in\Sigma$, we obtain that the time $t(x,z)$ becomes a function of the variable $x$ and the integral $\int_0^t(c_0x^2+c_2z^2)ds$ can be expressed as a function $q$ of the variable $x$. Observe that $q(x)>0$ for every $x$. Therefore we get that,

\begin{align}\label{General Formula y(t)for Neutral2}
y(t)&=y_0\ex{-\ell t}\ex{-\ell q(x)}.
\end{align}

Hence $y(t)$ decreases exponentially fast as before but with a faster rate. All together, we get that the map $\fun{D_1}{\Sigma}{S}$ is given by

\begin{align}
\label{Map N1 Neutral model2}
D_1(x,y,1) &=\on{N}^{\tau'_{\on{Neu}}(x)}(x,y,1) \nonumber \\
      &=\Big(1,y\ex{-\ell(\Abtwo{x}{\beta_2}+q(x))},\Ab{x}{\beta}),
\end{align}

\noindent where $\beta=\frac{\beta_0}{\beta_2}$, compare this with \eqref{Map N1 Neutral model1}. The statements form Observation~\ref{Obs map N1} for this new version of the map $D_1$ will also hold. To finish the construction of the Poincar\'e map we compose now with the map $P_2$. Therefore, the full return map $\fun{\Pneu}{\Sigma}{\Sigma}$ of $\Sigma$ is given by

\begin{equation}\label{P Neu map2}
\Pneu(x,y) =
\begin{cases}
\left(a\Ab{x}{\beta}-1,\;y\ex{-\ell(\Abtwo{x}{\beta_2}+q(x))}-\dfrac{1}{2}\right),        & \text{if} \; x \in (0, 1];\\[1em]

\left(1-a|\Ab{x}{\beta}|,\;y\ex{-\ell(\Abtwo{x}{\beta_2}+q(x))}+\dfrac{1}{2}\right),        & \text{if} \; x \in [-1, 0),\\
\end{cases}
\end{equation}

\noindent where $\beta=\frac{\beta_0}{\beta_2}\in(0,1)$.

The matrix $D\Pneu$ has eigenvalues $\lambda_1=a\beta |x|^{\beta-1}$ and $\lambda_2=\ex{-\ell(\Abtwo{x}{\beta_2}+q(x))}$. Since $\beta\in(0,1)$ we have that $\lambda_1\rightarrow\infty$ as $x\rightarrow0$. For the other eigenvalue we have that $\lambda_2<1$ and $\lambda_2\rightarrow0$ as $x \to 0$ since $(\Abtwo{x}{\beta_2}+q(x))\rightarrow\infty$ as $x\to0$. Thus the modified Poincar\'e map $\Pneu$ is hyperbolic when $x$ approaches the origin.

The properties stated before remain true for this new modified return map $\Pneu$ like the invariance of the stable foliation given by the vertical lines $\ell_v(c)$ under the map $\Pneu$. Hence, we can express again $\Pneu$ as $\Pneu(x,y)=(\fneu(x),\gneu(x,y))$, where $\fun{\fneu}{I\setminus\{0\}}{I}$ is again a Lorenz-like map with exponent $\beta$ (see \eqref{Lorenz map function}) and $\fun{\gneu}{(I\setminus\{0\})\times I}{I}$ satisfy the same properties stated in the previous section. The existence of a unique a.c.i.p and statistical properties such as exponential decay of correlations for observables with bounded variation for the Lorenz-like map $\fneu$ are well established, see for example \cite{Viana}.

\section{The stable foliation and the UNI condition}\label{Sec. stable foliation and UNI}

\subsection{Existence and regularity of the strong stable foliation}

In this subsection we will study the properties of the strong stable foliation $\mathcal{F}^{ss}$ for the neutral geometrical Lorenz model we built in Section~\ref{Sec. Poincare maps Neutral models}.

For the neutral geometrical Lorenz attractor, denoted by $\Lambda_{\on{Neu}}$, we consider the Lorenz attractor $\Lambda_{\on{Lor}}$ in an open neighbourhood $U$ of the origin. Instead of considering the linearised vector field we consider the vector field given by \eqref{Neutral system} with the parameters corresponding for model 1 and 2 described in Section~\ref{Sec. model 1} and~\ref{Sec. model 2}, respectively. More precisely, we take an open neighbourhood $U$ in which the cross-section $\Sigma$ is contained. Then the Dulac map from $\Sigma$ to $S$ has the form given by \eqref{Map N1 Neutral model1} and \eqref{Map N1 Neutral model2} for the models 1 and 2, respectively.

This modification yields a different flow time from the cross-section $\Sigma$ to $S$.
In the original Lorenz construction we have a logarithmic Poincaré return time but for these
modifications we have a polynomial Poincaré return time given by \eqref{pol flow time N1}. 
The rest of the construction, however, remains unchanged; that is, the flow
constructed from $S$ to $\Sigma$ is made by a composition of an expansion, a rotation 
and a translation. Therefore we have the same hitting time $\tau_2(x)$ and thus the full 
return time for the modified Poincar\'e map $\Pneu$ is given in
\eqref{EqReturnTimeNeu}. The modified Poincar\'e map $\fun{\Pneu}{\Sigma}{\Sigma}$ is
given by \eqref{P Neu Map1} and \eqref{P Neu map2} for the model 1 and 2, 
respectively. We saw that the lines in the $y$-direction, denoted by $\ell_v(c)$, in the cross-section 
$\Sigma$ form the stable foliation which is preserved by the return map $\Pneu$. Thus by 
quotienting out the stable direction we can rewrite the Poincar\'e map as a skew-product;
that is, $\Pneu(x,y)=(\fneu(x),\gneu(x,y))$, where $\fneu$ is a
one-dimensional Lorenz-like map.

\begin{lem}\label{lem:ev}
If $a_2b_0 < 9a_0b_2$,
then the eigenvalues of $DZ_{(x,y,z)}$ satisfy $0 < -\lambda_s < \lambda_u < -\lambda_{ss}$
for all $(x,y,z) \in U$.
\end{lem}

\begin{proof}
The derivative matrix $DZ_{(x,y,z)}$ of the vector field $Z$ is
$$
\begin{pmatrix}
        3a_0x^2+a_1y^2+a_0z^2 & 2a_1xy & 2a_2xz \\
        -2c_0\ell xy & -\ell(1+c_0x^2+c_2y^2) & -2c_2 \ell yz \\
         -2b_0xz & -2b_1yz & -(b_0x^2+b_1y^2+3b_2z^2)
       \end{pmatrix}.
$$
In finding the eigenvalues, we get $\lambda_{ss} = -\ell$ for Model 1 (i.e., $c_0=c_2=0$)
and $\lambda_{ss} = -\ell(1+c_0x^2+c_2z^2)$ for  Model 2 (i.e., $a_1=b_1=0$),
In both cases, the other two eigenvalues are
$\frac12\left(\tr(A) \pm \sqrt{\tr(A)^2 - 4 \det(A)}\right)$
for the submatrix
$$
A = \begin{pmatrix}
        3a_0x^2+a_1y^2+a_0z^2 & 2a_2xz \\
         -2b_0xz & -(b_0x^2+b_1y^2+3b_2z^2)
       \end{pmatrix}.
$$
To ensure that these eigenvalues are real and of opposite sign, it suffices to check that
$$
\lambda_u \lambda_s =
\det(A) = -(3a_0x^2+a_1y^2+a_0z^2)(b_0x^2+b_1y^2+3b_2z^2)+4a_2b_0x^2z^2 < 0.
$$
The worst case is when $y = 0$, so we consider the terms not including $y$ (and divide by $3$ for simplicity):
$$
(a_2b_0-3a_0b_2)x^2z^2 - a_0b_0 x^4 - a_2b_2 z^2 < 0.
$$
Divide by $z^4$ and introduce the new coordinate $u = x^2/z^2$:
\begin{equation}\label{eq:PQR}
-Pu^2 + Qu  -R < 0 \quad \text{ for }  \quad P = a_0b_0, Q = a_2b_0-3a_0b_2, R = a_2b_2.
\end{equation}
The left hand side is indeed negative for $u = 0$, and it is negative for all $u$ if the
equation $-Pu^2 + Qu  -R = 0$ has no real solution, i.e., if the discriminant is negative:
$$
0 > Q^2 - 4PR = (a_2b_0-3a_0b_2)^2 - 4a_0b_2a_2b_0 = 9a_0^2b_2^2 + a_2^2b_0^2 - 10 a_0b_2a_2b_0.
$$
Divide by $a_2b_0$ and use the new coordinate $\gamma = \frac{a_0b_2}{a_2b_0}$. Then we get the inequality
$$
0 > 9\gamma^2 - 10\gamma + 1 = (3\gamma - \frac{5}{3})^2 - \frac{16}{9},
$$
which is equivalent to $|\gamma - \frac{5}{9}| < \frac{4}{9}$.
That is, it fails if $\gamma \leq \frac19$ or $\gamma \geq 1$.
Now consider equality in\eqref{eq:PQR} and we divide by $a_2b_0$, which brings it to the form
$$
-\frac{a_0}{a_2} u^2 + (1-3\gamma) u - \frac{b_2}{b_0} = 0,
$$
with solutions
$$
u = \frac{a_2}{2a_0} \left( 1-3\gamma \mp \sqrt{ (1-3\gamma)^2 - 4\gamma} \right)
= \frac{a_2}{2a_0} \left( 1-3\gamma \mp \sqrt{ 9 \gamma^2 - 10 \gamma +1} \right).
$$
If $\gamma \geq 1$, then these solutions are negative, and since $u = x^2/z^2$, this means that there are no solutions $(x,y,z) \in U$. The remaining case $\gamma \leq \frac19$ is exactly the excluded case in the lemma.
This concludes the proof.
\end{proof}

Using our assumption $a_2b_0 < 9a_0b_2$ and Lemma~\ref{lem:ev}, we obtain
that the origin is the only point where we have $\lambda_{ss}=-\ell$ and $\lambda_s=\lambda_u=0$. Before continuing, we recall the definitions of a partially hyperbolic set and strongly dissipativity.

\bdn\label{DefPartiallyHyp}
Let $\Lambda$ be a compact invariant set of $X\in\vf{X}{r}(M)$, $c>0$ and $0<\lambda<1$. We say that $\Lambda$ has a $(c,\lambda)-$\textbf{dominated splitting} if the tangent bundle $T_{\Lambda}M$ has a $DX^t$-invariant splitting of sub-bundles
$$
T_{\Lambda}M=E^1\oplus E^2,
$$ 
\noindent such that for all $t>0$ and $x\in\Lambda$, we have

\begin{equation}\label{DominationSplitting}
\norm{DX^t|_{E_x^1}}{}\cdot\norm{DX^{-t}|_{E^2_{X^t(x)}}}{}<c\cdot\lambda^t.
\end{equation}
We say that $\Lambda$ is \textbf{partially hyperbolic} if it has a $(c,\lambda)-$dominated splitting such that $E^1$ is uniformly contracting; that is, for some $c>0$ and all $t>0$ and every $x\in\Lambda$ it holds 

\begin{equation}\label{ContractionSplitting}
\norm{DX^t|_{E_x^1}}{}<c\cdot\lambda^t.
\end{equation}

\noindent In this case we will denote $E^1$ by $E^s$ and call it the contracting direction. Also $E^2$ will be denoted by $E^{cu}$ and called the center-unstable direction.
\edn

\bdn\label{DissipativityFlow}
Let $G$ be a $\Cdiff{\infty}$ vector field on $\R^3$ with a Lorenz-like equilibrium; that is, the eigenvalues of $DG_p$ are real and satisfy $0<-\lambda_s<\lambda_u<-\lambda_{ss}$. We say that $G$ is \textbf{strongly dissipative} if the divergence of the vector field $G$ is strictly negative; that is, there is $c>0$ such that $\Div{G}(x)\leq-c$ for every $x$ and the eigenvalues of the singularity at $p$ satisfy the constraint 

\begin{equation}\label{ConstraintEigenvalues}
\lambda_u+\lambda_{ss}<\lambda_s.
\end{equation}
\edn

Figure~\ref{fig:Map N1} shows how the flows given by the Neutral model 1 and 2 send the lines in the $y$-direction in $\Sigma$ to lines in the $y$-direction in $S$. Thus, its derivative $D\on{N}^t$ preserves the lines in the $y$-direction. Furthermore, by the way the flow from $S$ to $\Sigma$ was constructed (Figure~\ref{fig:Lorenz}) we notice that horizontal lines in $S$; that is, parallel to the $y$-axis, are taken to parallel lines to the same axis in $\Sigma$. In other words, the flow from $S$ to $\Sigma$ preserves parallel lines to the $y$-axis. Since this flow is a composition of
a rotation, an expansion and a translation, the derivative of the flow also preserves planes orthogonal to the $y$-axis. From this we can deduce that the splitting $\R^3=E\oplus F$, where $E=\{0\}\times\R\times\{0\}$ and $F=\R\times\{0\}\times\R$, is preserved by the flow; \textit{i.e.}, $D\on{N}^t(E)=E$ and $D\on{N}^t(F)=F$ for any $t$, where $\on{N}^t$ is the flow obtained from Equation~\eqref{Neutral system} with the corresponding parameters for model 1 and 2. Since we have uniform contraction along $E$ ($\norm{D\on{N}^t|_{E_x}}{}\leq\ex{\lambda_{ss}t}$ with $\lambda_{ss}=-\ell<0$ for every $x\in U$) and domination of the splitting ($\norm{D\on{N}^t|_{E_x}}{}\cdot\norm{D\on{N}^{-t}|_{F_{\on{N}^t(x)}}}{}\leq\ex{\lambda_{ss}-\lambda_s t}$ with $\lambda_{ss}-\lambda_s<0$ for every $x\in U$) we 
can conclude that $U$ is partially hyperbolic. It is worth noticing 
that the origin is the only point that spoils the singular 
hyperbolicity condition (a set $A$ is singular hyperbolic if all its 
singularities are hyperbolic and it has volume expanding central
direction) since $J^{cu}_t(\bar{0})=|\det D\on{N}^t|_{F_0}|=\ex{(\lambda_u+\lambda_s)t}=1$ (recall $\lambda_s=0=\lambda_u$); that is, there is no area expansion along the subbundle $F$. Hence, $\Lambda_{\on{Neu}}$ is a partially hyperbolic attractor.

Theorem 6 in \cite{DecayCorrelationsContractingFibers} provides us with local strong-stable and center-unstable laminations $W^{ss}_\epsilon(x)$ and $W^{cu}_\epsilon(x)$, respectively, through the points $x\in\Lambda_{\on{Neu}}\setminus\{\bar{0}\}$. We note that both, $W^{ss}_\epsilon(x)$ and $W^{cu}_\epsilon(x)$, are embedded discs and hence sub-manifolds of $M$.
Also  $W^{ss}_\epsilon(x)$ is uniquely determined since $E^s$ is uniformly contracting. Corollary 6 in \cite{DecayCorrelationsContractingFibers} shows us that the local strong-stable lamination can be extended to an invariant foliation $\mathcal{F}^{ss}(x)$ of a open neighbourhood of $\Lambda_{\on{Neu}}$ with $\Cdiff{2}$ leaves and whose foliated charts are $\Cdiff{1}$. Moreover, the leaves are uniformly contracted by the action of the flow.


We note that $\Sigma$ is a $\Cdiff{2}$ embedded compact disk transversal to the flow $\on{N}^t$. Furthermore, $\Sigma$ is contained in the open neighbourhood $V$ of $\Lambda_{\on{Neu}}$. By Theorem 6 and Corollary 6 from \cite{DecayCorrelationsContractingFibers} we know that local strong-stable lamination $W^{ss}_\epsilon(x)$ extends to an invariant foliation $\mathcal{F}^{ss}(x)$. In this way, for $x\in\Sigma$ we define $W^{ss}(x,\Sigma)$ to be the connected component of $\mathcal{F}^{sc}(x)\cap\Sigma$, where $\mathcal{F}^{sc}(x)=\bigcup_{t\in\R}\on{N}^t(\mathcal{F}^{ss}(x))$ is the center-stable leaf. Since the flow $(\on{N}^t)_{t\in\R}$ is $\Cdiff{2}$, $W^{ss}(x,\Sigma)$ is a $\Cdiff{2}$ one-dimensional embedded curve for every $x\in\Sigma$ and their leaves form a $\Cdiff{1}$ foliation $\mathcal{F}^{ss}_\Sigma$ of $\Sigma$.

\begin{figure}[H]
\centering
\includegraphics[scale=.4]{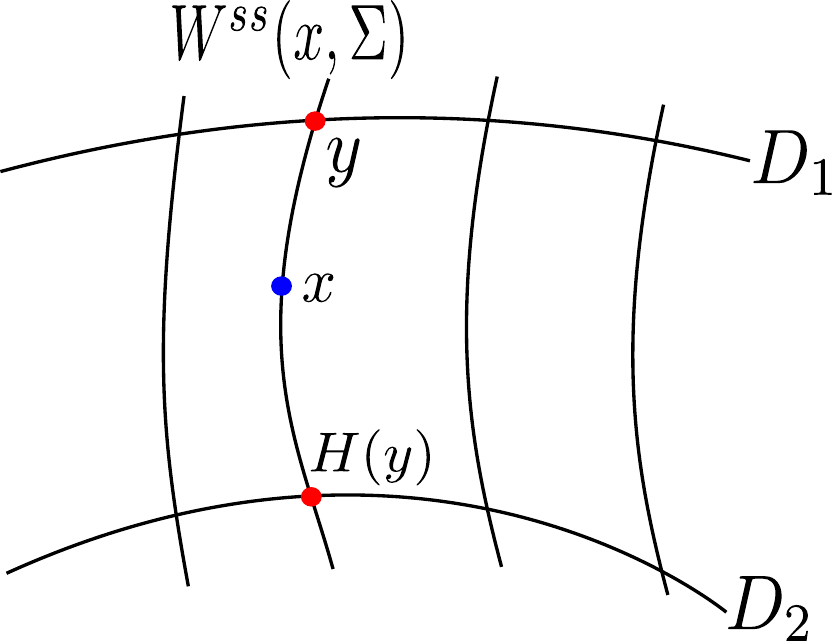}
\caption{The holonomy map.}
\label{holonomy map fig}
\end{figure}

Given a pair of embedded disks $D_1$ and $D_2$ in $\Sigma$ intersecting transversally a set $\{W^{ss}(x,\Sigma)\}_{x\in\Sigma}$ of stable leafs, the \textbf{holonomy map} $\fun{H}{D_1\cap W^{ss}(x,\Sigma)}{D_2\cap W^{ss}(x,\Sigma)}$ assigns to $y\in D_1\cap W^{ss}(x,\Sigma)$ the unique point in $h(y)\in D_2\cap W^{ss}(x,\Sigma)$, see Figure~\ref{holonomy map fig}.

From the developments in the partially hyperbolicity theory by Brin-Pesin \cite{BrinPesinPartially} and Pugh-Shub \cite{PughShubErgodicity} we have that the projection along leaves, also known as holonomies, between pair of transversal surfaces to $\mathcal{F}^{ss}$ have a Hölder continuous Jacobian with respect to Lebesgue measure. This implies a similar conclusion for the holonomies transversal to $\mathcal{F}^{cs}$. It follows that the holonomy between pairs of transversal curves to $\mathcal{F}_\Sigma^{ss}$ along the lines of $\mathcal{F}^{ss}_\Sigma$ can be thought of as interval maps having a Hölder Jacobian. Hence these holonomies are $\Cdiff{1+\alpha}$ for some $0<\alpha<1$. In this setting, the leaves $W^ss(x,\Sigma)$, with $x\in\Sigma$, determine a foliation $\mathcal{F}^{ss}_\Sigma$ of $\Sigma$ with transversal smoothness $\Cdiff{1+\alpha}$.

Therefore, we can assume that $\Sigma$ is the image of the unit square $I\times I$ under the action of a $\Cdiff{1+\alpha}$ diffeomorphism $h$, for some $0<\alpha<1$. Furthermore, $h$ sends vertical lines inside the leaves of $\mathcal{F}^{ss}_\Sigma$. The next step is to prove that the strong stable foliation $\mathcal{F}^{ss}(x)$ is not only $\Cdiff{1}$ but also $\Cdiff{1+\alpha}$, for some $0<\alpha<1$. This was done by Ara\'ujo, Melbourne and Varandas in \cite{RapidMixingLorenz}, stated as Lemma 2.2. This result is a consequence of domination of the splitting (Equation~\eqref{DominationSplitting}), uniform contraction along the stable direction (Equation~\eqref{ContractionSplitting}) and strong dissipativity (Definition~\ref{DissipativityFlow}). Therefore, we can conclude that the neutral geometrical flow has a strong-stable foliation $\mathcal{F}^{ss}(x)$ which is $\Cdiff{1+\alpha}$. Furthermore, the modified return map also has a strong-stable foliation $\mathcal{F}^{ss}_\Sigma$, whose transversal smoothness is $\Cdiff{1+\alpha}$, with $0<\alpha<1$.

A final remark concerning this subsection, Theorem 6 in \cite{DecayCorrelationsContractingFibers} is stated for a singular hyperbolic attractor. However, the conclusions and arguments still hold true if we consider a compact partially hyperbolic invariant set instead of a singular hyperbolic set. We will also like to mention that the situation for the origin is slightly different, since the splitting of the tangent bundle is given by a one-dimensional strong-stable direction $E^s$ and a two-dimensional center direction $E^c$. However, we are only concerned with the existence of the strong-stable manifold $W^{ss}(\bar{0})$, which follows from the theory of partial hyperbolicity.

\subsection{The UNI condition}\label{Sec. UNI condition}

The main goal of this subsection is to show that the stable and unstable manifolds of the modified geometrical model are jointly nonintegrable. The joint nonintegrability of stable and unstable foliations can be interpreted as follows: The stable and unstable foliation of an Anosov flow are always transversal, therefore, if they are jointly integrable, this provides us with a codimension one invariant foliation which is transversal to the flow direction. In contrast, if there exists a codimension one submanifold transversal to the flow direction, then this foliation must be subfoliated by both the stable and unstable foliations. Thus they must be jointly integrable. In this situation it is known \cite[Proposition 3.3]{FieldMelbourne} that the flow is semiconjugate to a suspension with a locally constant roof-function over a subshift of finite type. Such a flow need not mix! From the work of Ara\'ujo, Butterley and Varandas \cite{AxiomA}, we know that the joint nonintegrability of the stable and unstable manifolds is equivalent to the uniformly nonintegrability (UNI) condition. As we will see later in this section, the UNI condition, stated formally in Definition~\ref{UNI definition}, will ensure that the roof function of the suspension flow, which we will use to model the Neutral geometrical Lorenz flow, is not cohomologous to a constant function; that is, the time it takes for a point in the base to reach the roof of the suspension flow is not the same for every point. This property will guarantee the mixing properties.

In this chapter we will aim to prove the UNI condition. In \cite{RobustExpDecay} and \cite{RapidMixingLorenz} Ara\'ujo \textit{et al.\ } prove the UNI condition by exploiting the properties obtained by using hyperbolic times. From now on, when we talk about the neutral geometrical model, we will refer to both models we constructed in Section~\ref{Sec. Poincare maps Neutral models}.

Let $I=[-1,1]$ and $\fun{\fneu}{I}{I}$ be the one-dimensional Lorenz-like map obtained from the neutral geometrical model. We notice that $\{0\}$ is a nondegenerate\footnote{Nondegenerate here is meant for the derivative of the one-dimensional map $\fneu$ and shouldn't be confused with the degeneracy at the singularity $0$ concerning the eigenvalues being zero.} exceptional set for $\fneu$. From Theorem  4.3 in \cite{RobustExpDecay} we know that there are $\tilde{X}$ a neighbourhood of the singularity $0$, a countable partition $\tilde{\mathcal{Q}}$ of $\tilde{X}$ Lebesgue modulo zero into subintervals, a function $\fun{\tau}{\tilde{X}}{\Z^+}$ constant on partition elements and the induced map $\fun{\tilde{F}=\fneu^\tau}{\tilde{Q}}{\tilde{X}}$ which is a $\Cdiff{2}$ uniformly expanding diffeomorphism with bounded distortion for any $\tilde{Q}\in\tilde{\mathcal{Q}}$. 

The motivation behind taking the inducing scheme is that we aim to 
extend the Gibbs-Markov map $\tilde{F}$ to a two-dimensional Gibbs-
Markov map $F$ and build the suspension flow $F^t$ over the map $F$ 
with roof function $R$ (see Equation \eqref{roof function suspension}). 
This will allow us to use the results given in 
\cite{PolynomialDecayBalint} to deduce the decay of correlations for 
the suspension flow $F^t$ and ultimately for the geometrical neutral 
Lorenz flow $\on{N}^t$. Furthermore, we want to use the properties 
provided by the hyperbolic times and the bounded distortion of the map 
$\tilde{F}$, 
to obtain the bound stated in Proposition \ref{bound roof 
function neu case} for the roof function $R$, which will help us 
prove the UNI condition for the geometrical neutral flow.
We make now some observations regarding the induced map $\tilde{F}$.

\bobs\label{induced markov map tilde F inequalities obs}
The map $\tilde{F}$ is obtained by inducing $\fneu$ on the interval $\tilde{X}$, and the inducing time is given by the sum of a hyperbolic time with a non-negative integer bounded by $N$, where $N$ is such that $\bigcup_{i=1}^N(\fneu^i)^{-1}(\{0\})$ is $2\delta$-dense in $\tilde{X}$. Furthermore, $\tilde{F}$ is a full branch Markov map onto $\tilde{X}$ since $0$ has dense preimages under $\fneu$. For more details we refer the reader to \cite{RobustExpDecay}.

In addition to the bounded distortion and uniform expansion of $\tilde{F}$, we have the following inequalities for $\fneu$ as a consequence of hyperbolic times. Given $\sigma\in(0,1)$ and $c>0$, there is a constant $b>0$ such that:
\begin{enumerate}
\item[1.-] (Backward contraction) Let $\tilde{Q}(x)$ denotes the element of the partition $\tilde{\mathcal{Q}}$ containing $x$, for $y\in\tilde{Q}(x)$
\begin{align*}
|\fneu^i(y)-\fneu^i(x)|\leq b\sigma^{\tau(x)-i}|\tilde{F}(y)-\tilde{F}(x)|, \;\text{  }\; i=0,\ldots,\tau(x)-1.
\end{align*}

\item[2.-] (Slow recurrence to the singular point)
\begin{align*}
|\fneu^i(x)|\geq\sigma^{c(\tau(x)-i)}, \;\text{  }\; i=0,\ldots,\tau(y)-1.
\end{align*}
\end{enumerate}
\eobs

Following \cite{RapidMixingLorenz}, our next step is to construct a piecewise uniformly hyperbolic map $F$ with infinitely many branches, which covers $\tilde{F}$. First, let $W^{ss}_{\Pneu}(x)$ denote the stable leaf under the Poincar\'e map $\Pneu$ containing the point $x$, $\fun{\pi}{\Sigma}{I}$ the projection map. We define $X=\bigcup\{W^{ss}_{\Pneu}(x)\;|\;x\in\tilde{X}\}$ as the union of stable leaves along $\tilde{X}$.
We also extend the induced time $\tau$ to a function on $X$ denoted also by $\tau$ and given by $\tau(x)=\tau(\pi(x))$. We make the following important observation on the tails of $\tau$.

\bobs\label{exponential tails of m_X and m_tildeX obs}
The tails of the return time $\tau$ and its extension also denoted by $\tau$ are exponential (see \cite{RychlikBV}); \textit{i.e.}, there exists a constant $c>0$ such that $\mu_{\tilde{X}}(\tau>n)=\mathcal{O}(\ex{-cn})$ and $\mu_{X}(\tau>n)=\mathcal{O}(\ex{-cn})$, where $\mu_{\tilde{X}}$ and $\mu_{X}$ are the SRB-measures of the Gibbs-Markov maps $F$ and $\tilde{F}$, see below.
\eobs

Now, we construct $\fun{F}{X}{X}$ the Poincar\'e map by setting $F(x)=\Pneu^{\tau(\pi(x))}(x)$ for $x\in X$. Furthermore, let $\mathcal{Q}$ be the measurable partition of $X$ by taking $\bigcup\{W^{ss}_{\Pneu}(x)\;|\;x\in\tilde{Q}\}$ as its elements, with $\tilde{Q}\in\tilde{\mathcal{Q}}$. We will make use of $X$ and $F$ when making the model of the neutral geometrical Lorenz flow by a suspension flow.

It is standard \cite{RychlikBV} that the map $\tilde{F}$ has a unique a.c.i.p measure $\mu_{\tilde{X}}$ on $\tilde{X}$. Furthermore, $r_{\on{Neu}}\in\lp{1}{\mu_{\tilde{X}}}$ and there exists a unique invariant measure $\mu_X$ for $F$, $\mu_{\Sigma}$ for $\Pneu$ and $\mu_{I}$ for $\fneu$ satisfying $\pi_*(\mu_X)=\mu_{\tilde{X}}$, $\pi_*(\mu_\Sigma)=\mu_{I}$,
and also
\begin{eqnarray*}
\mu_\Sigma &=& \sum_{n\geq1}\sum_{i=0}^{n-1}(\Pneu^i)_*(\mu_X|\{\tau\circ\pi=n\}), \\
\mu_I &=& \sum_{n\geq1}\sum_{i=0}^{n-1}(\fneu^i)_*(\mu_{\tilde{X}}|\{\tau\circ\pi=n\}).
\end{eqnarray*}
Moreover, $\mu_X\ll\mu_\Sigma$ and $\mu_X(X)=1$, thus $\mu_\Sigma(X)>0$. We take the induced roof function $\fun{R}{X}{\R^+}$ given by

\begin{equation}\label{roof function suspension}
R(x)=\sum_{k=0}^{\tau(x)-1}\rneu{\Pneu^k(x)}.
\end{equation} 

Notice that $R$ is constant along stable leaves because $r_{\on{neu}}$ is constants along stable leaves. We also call $R$ the quotient induced roof function $\fun{R}{\tilde{X}}{\R^+}$. With this in mind we can state the definition of the UNI condition. First, we give the definition of suspension flow.

\bdn \label{SuspensionFlow} 
Let $(\Sigma,\nu)$ be a probability space and $\fun{P}{\Sigma}{\Sigma}$ an ergodic measure-preserving transformation. Let $\fun{r}{\Sigma}{\R^+}$ be a measurable (Hölder continuous) roof function. We define the \textbf{suspension space}\index{Space!Suspension} as $\Sigma^r=\{(x,u)\in \Sigma\times[0,r(x)]\}/\sim$, where $(x,r(x))\sim(P(x),0)$.
The \textbf{suspension flow}\index{Flow!suspension} $\fun{f_t}{\Sigma^r}{\Sigma^r}$ is given by $f_t(x,u)=(x,u+t)$ computed modulo identifications and the measure $\mu=\nu\times \lambda$, where $\lambda$ is the Lebesgue measure, is ergodic and $f_t$-invariant. In the finite measure case, we normalise by $\bar{r}=\int_\Sigma rd\mu$ so that $\mu=\frac{\nu\times \lambda}{\bar{r}}$ is a probability measure.
\edn

\bdn \label{UNI definition}
Let $\fun{R}{X}{\R^+}$ be a roof function as above, $\fun{F^t}{X^R}{X^R}$ the suspension flow built over $\fun{F}{X}{X}$, $R_n(x)=\sum_{i=0}^{n-1}R\circ F^i(x)$. Define $\fun{\psi_{h_1,h_2}=R_n\circ h_1-R_n\circ h_2}{X}{\R}$, for $h_1,\;h_2\in\mathcal{H}_n$; that is, inverse branches of $F^n$. Then the \textbf{UNI condition} holds if there exist $D>0$ and $h_1,\;h_2\in\mathcal{H}_{n_0}$ for some sufficiently large $n_0\geq1$, such that $|\psi'_{h_1,h_2}|\geq D$.
\edn

We saw in the previous subsection that the Poincar\'e map $\fun{\Pneu}{\Sigma}{\Sigma}$ has a strong stable foliation. We observe that the leaves of this foliation cross $\Sigma$, hence the induced map $F$ has a strong stable manifold $W^{ss}_F(x)=W^{ss}_{\Pneu}(x)$ that crosses $X$.
Ara\'ujo \textit{et al.}\ in \cite[Proposition 2.4]{RapidMixingLorenz} 
provide us with local unstable manifolds of uniform size for $F$ and 
defined almost everywhere, and by \cite[Proposition 2.4]{RapidMixingLorenz}
we obtain a local product structure, see Figure~\ref{fig:local product 
structure F}.

\begin{center}\vspace{1cm}
\includegraphics[width=0.45\linewidth]{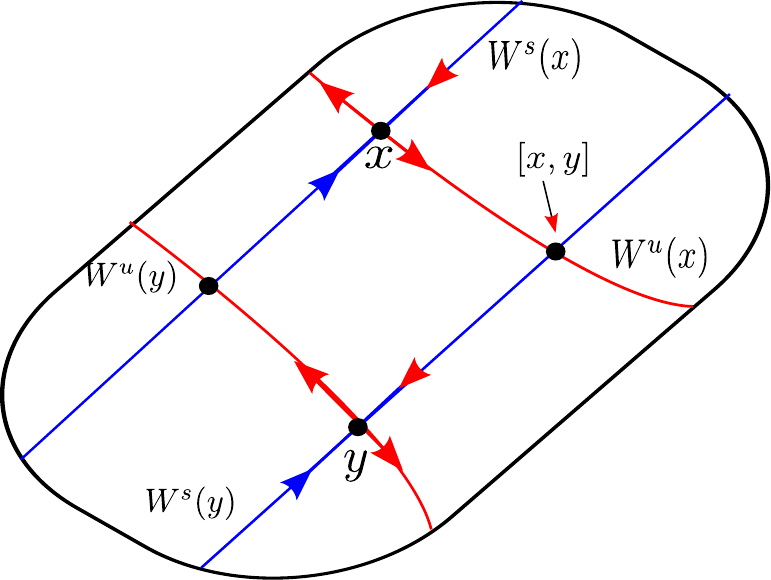}
\captionof{figure}{Local product structure for $F$.}
\label{fig:local product structure F}
\end{center}\vspace{1cm}

To use the arguments given in \cite{RapidMixingLorenz} for the temporal 
distortion function we need to adapt the proofs of the uniform bound of 
the derivative of the induced roof function $R$ and a bound on the flow 
time $r_{\on{Neu}}$. An adjustment is required because we are changing 
the flow time from $\Sigma$ to $S$. 
Recall that the original flow time for 
the geometrical Lorenz model is given by \eqref{EqReturnTimeLor}, 
whereas the modified flow time is given by \eqref{EqReturnTimeNeu}; 
that is, we change from logarithmic to polynomial. More precisely, we have the following propositions.

\bpro\label{uniforn bound derivative neu}
Let $\tilde{Q}\in\tilde{\mathcal{Q}}$ and $\fun{\tilde{F}}{\tilde{Q}}{\tilde{X}}$ as above. Denote by $\mathcal{H}$ the set of all inverse branches of $\tilde{F}$. Then we have
$$
\sup_{h\in\mathcal{H}}\sup_{x\in\tilde{X}}|D(R\circ h)(x)|<\infty.
$$
\epro 

\bdem
Let $\tilde{Q}\in\tilde{\mathcal{Q}}$ and $h\in\mathcal{H}$; that is, $\fun{h}{\tilde{X}}{\tilde{Q}}$ be an inverse branch of $\tilde{F}$ with inducing time $\tau=\tau(\tilde{Q})\geq1$ and fix $x\in\tilde{Q}$. We first observe that 
\begin{eqnarray*}
|D(R\circ h)(x)|&=& |DR(h(x))|\cdot|Dh(x)|=\frac{|DR(h(x))|}{|D\tilde{F}(h(x))|}\\&=&\left|\sum_{i=0}^{\tau-1}\frac{(Dr_{\on{Neu}}\circ\fneu^i)\cdot D\fneu^i}{D\tilde{F}}\circ h(x)\right|.
\end{eqnarray*}

From the construction of the inducing partition using hyperbolic times (for more details on hyperbolic times see  \cite{AlvesHypTimes} and \cite{RobustExpDecay}), we have that backward contraction and slow recurrence to the singular point, see Observation~\ref{induced markov map tilde F inequalities obs}; that is, there are constants $\sigma\in(0,1)$, $b_0>0$ and $c \in (0,\frac{1}{2}]$ such that
\begin{enumerate}
\item[1.-] (Backward contraction) For $y\in\tilde{Q}(x)$
\begin{align*}
|\fneu^i(y)-\fneu^i(x)|\leq b_0\sigma^{\tau-i}|\tilde{F}(y)-\tilde{F}(x)|, \;\text{  }\; i=0,\ldots,\tau-1.
\end{align*}

\item[2.-] (Slow recurrence to the singular point)
\begin{align*}
|\fneu^i(x)|\geq\sigma^{c(\tau-i)}, \;\text{  }\; i=0,\ldots,\tau-1.
\end{align*}
\end{enumerate}

Notice that

\begin{equation}\label{backward contraction bound for fneu}
\left|\frac{D\fneu^i}{D\tilde{F}}\circ h(x)\right|\leq b_0\sigma^{\tau-i}\;\text{  }\; i=0,\ldots,\tau-1,
\end{equation}

\noindent where the inequality follows from the backward contraction.

Moreover, from the slow recurrence to the singularity we also get the following inequality;

\begin{equation}\label{slow recurrence bound for fneu}
|(Dr_{\on{Neu}}\circ\fneu^i)\circ h(x)|\leq b_1\sigma^{-c(1+\frac{1}{\beta_2})(\tau-i)}\;\text{  }\; i=0,\ldots,\tau-1,
\end{equation}

\noindent for some constant $b_1>0$. Altogether, Equations \eqref{backward contraction bound for fneu} and \eqref{slow recurrence bound for fneu} imply that

\begin{equation}
|D(R\circ h)(x)|\leq b\sum_{i=0}^{\tau-1}\sigma^{si},
\end{equation}
\noindent where $s=1-c(1+\frac{1}{\beta_2})$. Since $0<c\leq\frac{1}{2}$ and $\beta_2>2$, we have $s<1$. Therefore the sum converges, and we have that
$$
\sup_{h\in\mathcal{H}}\sup_{x\in\tilde{X}}|D(R\circ h)(x)|<\infty.
$$
\edem

For $x,\;y\in\tilde{X}$ we define the separation time $s(x,y)$ as the least integer $n\geq0$ such that $\tilde{F}^n(x)$ and $\tilde{F}^n(y)$ are in different partition elements of $\mathcal{Q}_0$. For given $0<\eta<1$, the symbolic metric is defined on $\tilde{X}$ as $d_\eta(x,y)=\eta^{s(x,y)}$. Finally, we set $r_{\on{Neu}}^{(k)}(x)=\sum_{i=0}^{k-1}\rneu{\Pneu^i(x)}$.

\bpro\label{bound roof function neu case}
There exists $B>0$ such that for all $x,\; y\in\tilde{X}$ with $s(x,y)\geq1$ and $0\leq k\leq\tau(x)=\tau(y)$ we have $|r^{(k)}_{\on{Neu}}(x)-r^{(k)}_{\on{Neu}}(y)|\leq B|\tilde{F}(x)-\tilde{F}(y)|^\epsilon$. Consequently, there is $\eta\in(0,1)$ such that $|R|_\eta<\infty$, where $|R|_\eta=\sup_{x\neq y}\frac{|R(x)-R(y)|}{d_\eta(x,y)}$ denotes the Lipschitz constant of the quotient induced roof function $\fun{R}{\tilde{X}}{\R^+}$ with respect to $d_\eta$. Moreover, $|\tilde{F}(x)-\tilde{F}(y)|\leq Bd_\eta(x,y)$.
\epro

\bdem
For convenience, in this proof $f$ will denote $\fneu$. Let $x,\;y\in\tilde{X}$ such 
that $s(x,y)=n\geq1$ and $0\leq k\leq\tau(x)=\tau(y)$. Thus $y\in\tilde{Q}^n(x)$, 
where $\tilde{Q}^n(x)=\displaystyle\bigvee_{i=0}^{n-1}(\tilde{F}^i)^{-1}(\tilde{Q}(x))$ is the $n$th refinement of $\tilde{Q}(x)$, and so $\tau(\tilde{F}^i(x))=\tau(\tilde{F}^i(y))$ for $i=0,\ldots,n-1$. Hence, the choice of the cross-
section assures that $r_{\on{Neu}}$ is constant along stable leaves and that $\rneu{x}=|x|^{-\frac{1}{\beta_2}}h(x)+\tau_2(x)$, where $h(x)$ is bounded and bounded
away from zero.
In fact, $h(x)$ is of the form $h_0 +  \mathcal{O}(|x|^\gamma)$ where $h_0$ is a
positive constant and $\gamma > 0$ depending on whether the higher order terms are 
consider in \eqref{FlowHor} or not. Also $h(x)$ is differentiable for $x > 0$, 
because the Dulac map is differentiable. Then we can write

\begin{eqnarray*}
\label{Bound deriv. inequ 1}
|R(x)-R(y)| &\leq& \sum_{i=0}^{\tau(x)-1}|r_{\on{Neu}}(f^i(x))-r_{\on{Neu}}(f^i(y))|\nonumber\\ &\leq& \sum_{i=0}^{\tau(x)-1}\left||f^i(x)|^{-\frac{1}{\beta_2}}h(f^i(x))-|f^i(y)|^{-\frac{1}{\beta_2}}h(f^i(y))\right| \\&&+|\tau_2(f^i(x)-\tau_2(f^i(y))|.
\end{eqnarray*}

We first notice that,

\begin{eqnarray}
\label{bound tau}
\left|\tau_2(f^i(y))-\tau_2(f^i(x))\right|
  &\leq& \norm{\tau_2}{\epsilon}\left|f^i(y)-f^i(x)\right|^\epsilon \nonumber\\
  &\leq& \sigma^{\epsilon(\tau(x)-i)}\norm{\tau_2}{\epsilon}\left|\tilde{F}^i(y)-\tilde{F}^i(x)\right|^\epsilon.
\end{eqnarray}

The second inequality in \eqref{bound tau} follows from Observation~\ref{induced markov map tilde F inequalities obs}-$1$. Now, denote $\left| |f^i(x)|^{-\frac{1}{\beta_2}}h(f^i(x)) - |f^i(y)|^{-\frac{1}{\beta_2}}h(f^i(y))\right|$ by $A$, then we obtain that,

\begin{eqnarray}
\label{A bound}
 A &\leq& \left| |f^i(x)|^{-\frac{1}{\beta_2}} - |f^i(y)|^{-\frac{1}{\beta_2}}\right| |h(f^i(x))| \nonumber \\
  && + |f^i(y)|^{-\frac{1}{\beta_2}} |h(f^i(x))-h(f^i(y))|.
\end{eqnarray}

We notice that $\left| |f^i(x)|^{-\frac{1}{\beta_2}} - |f^i(y)|^{-\frac{1}{\beta_2}}\right|$ is bounded. Indeed we have the following:

\begin{eqnarray}
\label{bound 1}
  \left| |f^i(x)|^{-\frac{1}{\beta_2}} - |f^i(y)|^{-\frac{1}{\beta_2}}\right|
  &=& |f^i(x)|^{-\frac{1}{\beta_2}} \left| 1 - \left(\frac{ |f^i(y)|}{|f^i(x)|}
\right)^{-\frac{1}{\beta_2}}\right| \nonumber \\
  &\leq& |f^i(x)|^{-\frac{1}{\beta_2}} \left| 1 - \left( 1 + \frac{
|f^i(y)-f^i(x)|}{|f^i(x)|} \right)^{-\frac{1}{\beta_2}}\right| \nonumber \\
  &\leq&  |f^i(x)|^{-\frac{1}{\beta_2}} \frac{C_0}{\beta_2} \frac{
|f^i(y)-f^i(x)|}{|f^i(x)|} \nonumber \\
  &=& \frac{C_0}{\beta_2} |f^i(x)|^{-(1+\frac{1}{\beta_2})} |f^i(y)-f^i(x)| \nonumber \\
  &\leq& \frac{C_0}{\beta_2} \sigma^{(1-c(1+\frac{1}{\beta_2}))(\tau(x)-i)}\left|\tilde{F}(y)-\tilde{F}(x)\right|,
\end{eqnarray}
where $C_0>0$ and $0<c\leq\frac{2}{3}$. The first and second inequalities follow from the Bernoulli inequality and from Observation~\ref{induced markov map tilde F inequalities obs}-$1$ and $2$, respectively. Since $h(f^i(x))$ is bounded and by \eqref{bound 1}, we can bound the first term in the sum of \eqref{A bound} by $\frac{C_0}{\beta_2} \sigma^{(1-c(1+\frac{1}{\beta_2}))(\tau(x)-i)}\left|\tilde{F}(y)-\tilde{F}(x)\right|$.

To finish the proof it remains to find a bound for

\begin{equation}
\label{bound h terms}
|f^i(y)|^{-\frac{1}{\beta_2}} |h(f^i(x))-h(f^i(y))|. 
\end{equation}

Notice that $|h(f^i(x))-h(f^i(y))|\approx h'(\xi)|f^i(x)-f^i(y)|$, with $f^i(x)<\xi<f^i(y)$. Hence, by Observations~\ref{induced markov map tilde F inequalities obs}-$1$ and $2$, we have that \eqref{bound h terms} is bounded by $C_0h'(\xi)\sigma^{(1-\frac{c}{\beta_2})(\tau(x)-i)}|\tilde{F}(y)-\tilde{F}(x)|$. Assuming that $h'(\xi)$ is bounded and combining all the previous bounds we have that,

\begin{eqnarray*}
|R(x)-R(y)|
  &\leq& C\sum_{i=0}^{\tau(x)-1}[ (\sigma^{s(\tau(x)-i)}+\sigma^{u(\tau(x)-i)})|\tilde{F}(y)-\tilde{F}(x)|\\
  & & +\ \sigma^{\epsilon(\tau(x)-i)}\norm{\tau_2}{\epsilon}|\tilde{F}^i(y)-\tilde{F}^i(x)|^\epsilon]\\
  &\leq& B |\tilde{F}^i(y)-\tilde{F}^i(x)|^\epsilon,
\end{eqnarray*} 
for some constant $B>0$
where $s=1-c(1+\frac{1}{\beta_2}) > 0$ and $u=1-\frac{c}{\beta_2}$. As in the previous proof, the sum converges since $\beta_2>2$, $0<c\leq\frac{2}{3}$ and hence $0<s,u<1$. This establishes what we were aiming to prove.

One caveat we have to make here is that this argument is only valid if we assume boundedness for $h'$. This is true if the higher order terms in \eqref{FlowHor} are not present. If higher order terms are present, boundedness of $h'$ is plausible, but since the required perturbation argument in \cite{Regular} is less constructive so as to immediately derive this boundedness, we will try to convince the reader with the numeric analysis performed at the end of this work that this boundedness is indeed true. Thus, rather than a rigorous proof we give a combination of mathematical arguments and numerical verification.
\edem

Now let $Q\in\mathcal{Q}$ be a partition elements for $F$. The temporal distortion function $\fun{T}{Q\times Q}{\R}$ is defined almost everywhere by,

\begin{eqnarray}\label{temporal distortion fun T neu}
T(x,y)&=&\sum_{i=-\infty}^\infty[\rneu{\Pneu^i(x)}-\rneu{\Pneu^i([x,y])}-\rneu{\Pneu^i([y,x])}\nonumber\\
&& +\ \rneu{\Pneu^i}(y)]\nonumber\\
&=&\sum_{i=-\infty}^{-1}[\rneu{\Pneu^i(x)}-\rneu{\Pneu^i([x,y])}-\rneu{\Pneu^i([y,x])}\nonumber\\
&&+\ \rneu{\Pneu^i(y)}],
\end{eqnarray}
where $[x,y]$ is the local product of $x$ and $y$ (see Figure~\ref{fig:temp distortion fun}). The second inequality follows from the property of $r_{\on{Neu}}$ of being constant along stable leaves.

\begin{center}\vspace{1cm}
\includegraphics[width=0.55\linewidth]{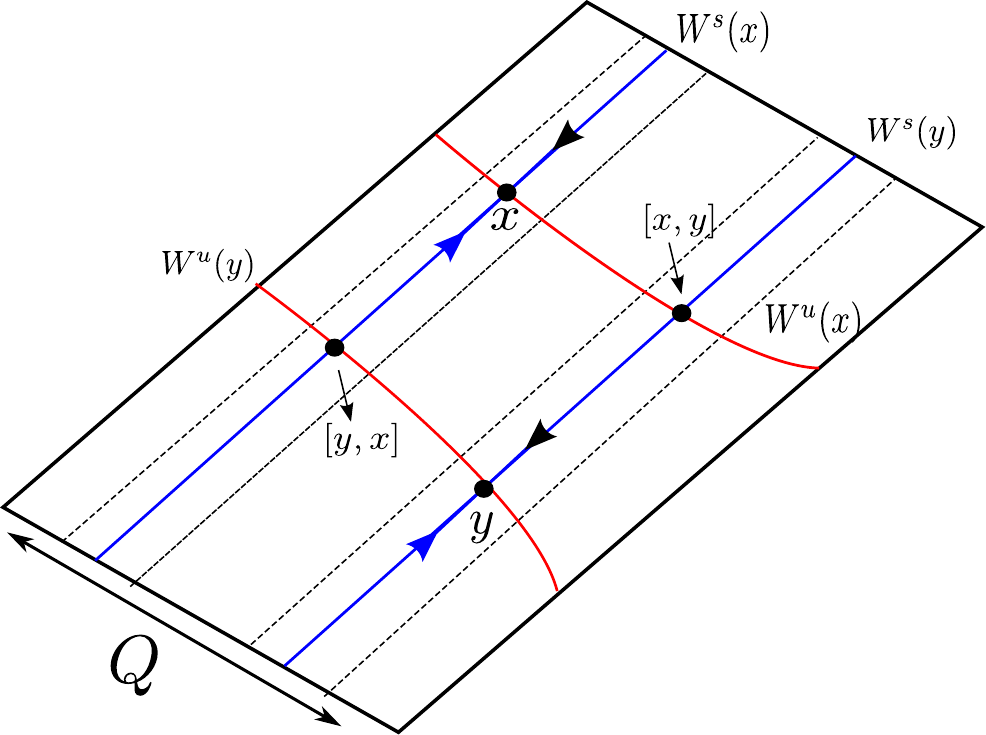}
\captionof{figure}{Local product structure for $F$.}
\label{fig:temp distortion fun}
\end{center}\vspace{1cm}

Now, for every $x,\; y\in X$ in the same unstable manifold for $\fun{F}{X}{X}$ we define

\begin{align}
T_0(x,y)=\sum_{i=1}^\infty[r_{\on{Neu}}(\Pneu^{-i}(x))-r_{\on{Neu}}(\Pneu^{-i}(y))].
\end{align}

The continuity and other properties of $T_0$ are stated in \cite[Lemma 3.1]{RapidMixingLorenz}. Furthermore, we can rewrite the temporal distortion function $T(x,y)$ in terms of $T_0$; that is,

$$
T(x,y)=T_0(x,[x,y])+T_0(y,[y,x]).
$$

The main result concerning the temporal distortion function establishes the joint nonintegrability of the stable and unstable foliations for the flow by proving that the temporal distortion function $T$ is not identically zero, that is, there is $Q\in\mathcal{Q}$ and $x,\; y\in Q$ such that $T(x,y)\neq0$ and is stated and proven in \cite[Theorem 3.4]{RapidMixingLorenz}.

We adjusted the proof of the uniform bound of the derivative of the induced roof function for the geometrical neutral Lorenz model since we want to apply \cite[Proposition 7.4]{ExpMixTeich} in order to use the same arguments given in \cite[Corollary 4.3]{ExpDecayNonUniform}. Thus, we get the UNI condition for the geometrical neutral Lorenz flow. For completeness we state it in the following theorem.

\bteo\label{UNI condition modified Lorenz}
The UNI condition holds for the geometrical neutral Lorenz flow.
\eteo

For fixed $x\in X$, we define the map $\fun{h}{W^{u}_F(x)}{\R}$ given by
$$
h(y)=T(x,y)=T_0(x,[x,y])+T_0(y,[y,x]),
$$
the map $h$ is $\Cdiff{1}$. Furthermore, there exists a nonempty open set $U\subset W^u_F(x)$ such that $h|_U$ is a $\Cdiff{1}$ diffeomorphism. For the proofs of the properties of the map $h$ see Proposition 3.6 and Corollary 4.7 in \cite{RapidMixingLorenz}.

The next result will be of great help in proving the decay of correlations for the geometrical neutral Lorenz flow. The original statement involves the geometrical Lorenz flow, but the same arguments can be used to prove the same result for our setting. Before stating the result we give a definition.

\bdn\label{finite subsystem}
Let $X$ and $F$ be as in the beginning of Section 6.3. A subset $Z_0\subset X$ is called a \textbf{finite subsystem}\index{finite subsystem} of $X$ if $Z_0=\bigcap_{n\geq0}F^{-n}Z$, where $Z$ is the union of finitely many partition elements of $X$.
\edn  

Let $Q_1$ and $Q_2\in\mathcal{Q}$ be two partition elements and consider $Q=Q_1\cup Q_2$. We define the finite subsystem $Q_0=\bigcap_{n\geq0}F^{-n}Q$, then we have the following:

\bpro\cite[Proposition 3.8]{RapidMixingLorenz}\label{temporal distortion function lower box dim positive}
$\;$ For the finite subsystem $Q_0$, the set $T(Q_0\times Q_0)$ has positive lower box dimension.
\epro

We will like to end this Section with the following remark. To establish their results on decay of correlations, B{\'a}lint et al.\ in \cite{PolynomialDecayBalint} and Melbourne in \cite{MelbournePol} assumed a very important, yet technical property namely, absence of approximate eigenfunctions. They also provide some criteria that guarantees the absence of approximate eigenfunctions. The first one, involves the temporal distortion function providing a nonintegrability condition. This criteria is given by Proposition \ref{temporal distortion function lower box dim positive}; that is, when the temporal distortion function is not identically zero. In other words, when the UNI condition is satisfied. The second one, is a Diophantine condition on the periods of three periodic solutions \cite{DolDio}, which is satisfied with probability one. It is important to remark that from these criteria, the UNI condition is robust while the Diophantine condition is not.

\section{Decay of correlations}\label{Sec. proof of decay of correlations Neutral}

In this section we prove Theorem~\ref{Poly decay Neutral geom Lorenz flow Thm} for the first two neutral models. The third and more general model will be analysed in Section~\ref{Sec. numerics}. We will use the results of B\'alint \textit{et al.\ } \cite{PolynomialDecayBalint} to prove our theorem. In \cite{PolynomialDecayBalint} polynomial decay of correlations for non-uniformly hyperbolic flows is proven under absence of approximate eigenfunctions. Let us start by giving the description of a non-uniformly hyperbolic flow described in \cite{PolynomialDecayBalint}.

First, we observe that the geometrical neutral Lorenz flow $\fun{\on{N}^t}{\Lambda_{\on{Neu}}}{\Lambda_{\on{Neu}}}$, where $\Lambda_{\on{Neu}}$ is the geometrical neutral Lorenz attractor, can be modelled as the suspension flow $\fun{S^t}{\Sigma^{r_{\on{Neu}}}}{\Sigma^{r_{\on{Neu}}}}$ over the Poincar\'e map $\Pneu$ with base the cross-section $\Sigma$ and roof function $r_{\on{Neu}}$ from \eqref{EqReturnTimeNeu}. However, in order to use the results of \cite{PolynomialDecayBalint}, we take the alternative model $\fun{F^t}{X^R}{X^R}$, where $X\subset \Sigma$ is a cross-section to the flow with nice hyperbolic structure (local product structure) and with induced roof function $\fun{R}{X}{\R^+}$ given by $R(x)=\sum_{k=0}^{\tau(x)-1}\rneu{\Pneu^k(x)}$, see Section~\ref{Sec. UNI condition}. Then the suspension flow $F^t$ built over the uniformly hyperbolic map $\fun{F}{X}{X}$ is identical to the suspension flow $S_t$, thus $F_t$ is an extension of the underlying flow, namely the neutral geometrical Lorenz flow. Within this framework, $\on{N}^t$ is called in \cite{PolynomialDecayBalint} a \textbf{non-uniformly hyperbolic flow}.

Under suitable conditions it can be shown that the suspension flow $F^t$ described above is a Gibbs-Markov flow \cite[Section 6]{PolynomialDecayBalint}. Therefore, the mixing rates for non-uniformly hyperbolic flows can be deduced from the corresponding results for Gibbs-Markov flows, see \cite[Corollary 8.1]{PolynomialDecayBalint}.

For observables $v$ and $w$, let $\rho_t(v,w)$ denote the decay of correlations of the geometrical neutral Lorenz flow; that is,
\begin{align}
\rho_t(v,w)=& \Bigg|\int v\cdot w\circ N^td\mu-\int vd\mu\int wd\mu\Bigg|,
\end{align}

\noindent where $\mu$ is the SRB measure of $\on{N}^t$. Before giving the proof of Theorem~\ref{Poly decay Neutral geom Lorenz flow Thm} for the neutral models 1 and 2, we will give the definitions of the space of observables.

Let $(M,d)$ be a metric space with $\on{diam}(M)\leq1$ and define a flow $\fun{T^t}{M}{M}$ on $M$. We fix $\eta\in(0,1]$ and for a given observable $\fun{v}{M}{\R}$ we define
$$
|v|_{C^\eta}=\sup_{x\neq y}\frac{|v(x)-v(y)|}{d(x,y)^\eta},
$$ 

\noindent and the norm $\norm{v}{{C^\eta}}=|v|_\infty+|v|_{C^\eta}$. We define the Banach space of Hölder functions on $M$ by $C^\eta(M)$; \textit{i.e.},
$$
C^\eta(M)=\{\fun{v}{M}{\R}\;|\;\norm{v}{{C^\eta}}<\infty\}.
$$

Furthermore, let 
$$
|v|_{C^{0,\eta}}=\sup_{\substack{x\in M\\t>0}}\frac{|v(T^t(x))-v(x)|}{t^\eta}
$$

\noindent and define $\norm{v}{{C^{0,\eta}}}=|v|_\infty+|v|_{C^{0,\eta}}$. We denote the space of Hölder observables in the flow direction by

$$
C^{0,\eta}(M)=\{\fun{v}{M}{\R}\;|\;\norm{v}{{C^{0,\eta}}}<\infty\}.
$$ 

We will say that an observable $\fun{w}{M}{\R}$ is \textbf{differentiable in the flow direction} if
$$
\partial_tw=\lim_{t\rightarrow0}\frac{w\circ T^t-w}{t}
$$
exists pointwise. Let $\norm{w}{{C^{m,\eta}}}=\sum_{k=0}^m\norm{\partial_t^kw}{{C^\eta}}$. We will denote the space of observables that are $m$-times differentiable in the flow direction by

$$
C^{m,\eta}(M)=\{\fun{w}{M}{\R}\;|\;\norm{w}{{C^{m,\eta}}}<\infty\}.
$$ 

For a Borel set  $X\subset M$ we define $C^\eta(X)$ as above by using the restriction of the metric $d$ to $X$.

\tbdem

We first note that the suspension flow $F^t$ projects to a quotient suspension semiflow $\fun{\tilde{F}^t}{\tilde{X}^R}{\tilde{X}^R}$, where $\tilde{X}$ is the quotient space obtained from $X$ after quotienting out the stable leaves and $\fun{\tilde{F}}{\tilde{X}}{\tilde{X}}$ is a Gibbs-Markov map. Proposition~\ref{bound roof function neu case} ensures that the following inequality holds.

\begin{align}\label{Gibbs-Markov skew roof function inequality}
|\varphi(x)-\varphi(y)|\leq C\gamma^{s(x,y)}\inf_{Q_i}\varphi && \text{for all }x,\; y\in\tilde{Q}_i,\; i\geq1,
\end{align}

where $\{\tilde{Q}_i\}_{i\geq0}$ is a countable Lebesgue modulo zero partition into subintervals, $0<\gamma<1$ and $s(x,y)$ is the separation time.

Therefore, we have that $\tilde{F}^t$ is a Gibbs-Markov semiflow and consequently that $F^t$ is a Gibbs-Markov flow. Then the conclusion follows from \cite[Corollary 8.1]{PolynomialDecayBalint}.

There are still four details concerning the hypothesis in \cite[Corollary 8.1]{PolynomialDecayBalint} that we have not mentioned yet. The first one is regarding condition (H); for us this condition holds automatically since $R$ is constant along stable leaves.

The second concerns the absence of approximate eigenfunctions for $F^t$. Melbourne gave in \cite[Chapter 5]{MelbournePol} sufficient conditions for the absence of approximate eigenfunctions, namely the existence of a finite subsystem with positive lower box dimension. Hence, it follows from Proposition~\ref{temporal distortion function lower box dim positive} and Lemma 8.9 in \cite{PolynomialDecayBalint} that $F^t$ has absence of approximate eigenfunctions.

The third concerns the tails of $R$; that is, we want to estimate $\mu_X(R>t)$. From Observation~\ref{exponential tails of m_X and m_tildeX obs} we know that $\mu_X(\tau>n)$ has exponential tails; \textit{i.e.}, there exists a constant $c_0>0$ such that $\mu_X(\tau>n)=\mathcal{O}(\ex{-c_0n})$. Moreover, by Theorem~\ref{Estiamete tails of return} we have that $\mu_X(r_{\on{Neu}}>t)$ has polynomial tails; that is, there is a constant $c_1>0$ such that $\mu_X(r_{\on{Neu}}>t)=\mathcal{O}(c_1t^{-\beta_2})$, where $\beta_2=\frac{a_2+b_2}{2b_2}$. Then by \cite[Proposition 5.1]{PolynomialDecay} we have that $\mu_X(R>t)=\mathcal{O}((\ln t)^{\beta_2}t^{-\beta_2})$.

The fourth and last detail concerns how to improve the estimates for $\mu_X(R>t)$ and remove the logarithmic term. For this, we make use of \cite[Lemma 4.1]{balint2010decay}. There the settings is made for infinite horizon planar periodic Lorentz gases, for that setting the tails of the flow time (in this work denoted by $r_{\on{Neu}}$) is of order $\mathcal{O}(t^{-2})$. By replacing in \cite[Lemma 4.1]{balint2010decay} the order of the tails from $\mathcal{O}(t^{-2})$ to $\mathcal{O}(c_1t^{-\beta_2})$ we can use the same proof to remove the logarithmic term. Hence, we have that $\mu_X(R>t)=\mathcal{O}(t^{-\beta_2})$. With this we conclude our proof.
\tedem

The natural questions is about the lower bounds, i.e., if the bounds
given in this theorem are sharp.

Although we definitely think they are, the currently available
literature is insufficient to conclude this, although the margin is fairly narrow.

In \cite{Gouezel}, the renewal operator methods are developed to get
such lower bounds, but his paper is for maps, not flows.
Melbourne and Terhesiu come the closest in \cite{MT17}, where they
consider suspension semiflows with polynomial roof functions over Gibbs-Markov base maps, and indeed, their results imply polynomial mixing for the flow $\fun{F^t}{X^R}{X^R}$, for H\"older observables that are constant on the stable fibers. The step from this suspension flow to the actual flow $\fun{\on{N}^t}{\Lambda_{\on{Neu}}}{\Lambda_{\on{Neu}}}$, however, is not trivial at all. (Here the effect of the second inducing step needs to be undone, in a way).

In \cite{PolynomialDecay} this step was taken for discrete suspensions (i.e., Young towers), specifically for billiard maps, but not for flows. However, the lower bounds that we obtain for $\fun{F^t}{X^R}{X^R}$ as a corollary of the results in \cite{MT17} are sufficient to prove stable laws (with exponent $1/\beta_2 \in (1,2]$) for the
neutral Lorenz flow $\fun{\on{N}^t}{\Lambda_{\on{Neu}}}{\Lambda_{\on{Neu}}}$, cf.\
\cite{BTT21, Volume}.

\section{Numerical analysis and results}\label{Sec. numerics}

In this section we provide the results of the numerical approximation we obtained for the exponent $\beta$ of the Dulac map (see Equation~\eqref{Dulac Map formula}) and the exponent $\beta_2$ of the tails of the return map (see Theorem~\ref{Estiamete tails of return}) for the two-dimensional setting (the setting of \cite{Regular}) and for the 3-dimensional setting concerning this work.

\subsection{Numerics of the $2$-dimensional case}
To verify the existing theoretical asymptotics from \cite{Regular}, we will start the numerical analysis in the $2$-dimensional case; that is, we consider the framework of \cite{Volume} and \cite{Regular}. There, the following neutral form was considered:

\begin{equation}
\label{Normal Flow 2D} 
\left( \begin{array}{ccc}
\dot{x}\\
\dot{y}\\
\end{array}\right)=
\left( \begin{array}{ccc}
x(a_0x^\kappa + a_2z^\kappa)\\
-y(b_0x^\kappa + b_2y^\kappa)\\
\end{array}\right)+\mathcal{O}(4),  
\end{equation}

where $a_0,\;a_2,\;b_0$ and $b_2>0$ and $\Delta:=a_2b_0-a_0b_2\neq0$. For simplicity we let $\kappa=2$. For the analysis of the Dulac map close to the neutral equilibrium of Equation~\eqref{Normal Flow 2D}, we take an unstable leaf $W^u(0,y_0)$ and a stable leaf $W^s(x_0,0)$, then the Dulac map $\fun{D}{W^u(0,y_0)}{W^s(x_0,0)}$, shown in Figure~\ref{fig:Dulac map}, assigns the firs intersection of the integral curve through $(x,y_0)$ with the stable leaf $W^s(x_0,0)$, where $x\in W^u(0,y_0)$ and $T$ is the flow time; that is, the exit time.

For the setting considered in this work, we will perform the numerical experiments with $x_0=1$. In order to corroborate the estimates of the Dulac map given by Equation~\eqref{Dulac Map formula}, we expect the numerical experiments to show us that

\begin{equation}
\beta\approx\frac{\ln(y)}{\ln(x)}.
\end{equation} 

We actually will see that 

\begin{equation}\label{2D Beta Adjust equation}
\beta=\frac{\ln(y)}{\ln(x)}-\frac{\ln(c(y_0))}{\ln(x)}+\mathcal{O}\Big(\frac{1}{\ln(x)}\Big).
\end{equation} 

From \cite{Volume} and \cite{Regular} we know that the constants $c(y_0)$ are given by a specific formula which is not easy to compute. For this reason, we decided to use the least-squares method to calculate these constants.

For the numerical experiments we will take different values of $\beta$ and  250 points $x\in[\num{1.0e-5},\num{1.0e-4}]$ at the unstable leave $W^u(0,y_0)$ with $y_0=1.0$. The integration method we will use for the numerical experiments concerning this work is the so-called Radau quadrature method, to deal with the numerical complications of integrating near a neutral stationary point, see \cite{Numerics}.

Figure~\ref{fig:2D beta approx} $a)$, $b)$ show us the approximation of $\beta$ (the red graph), the adjusted approximation of $\beta$ (green graph), and the theoretical value of $\beta$ (blue graph) for $\beta=0.266$ and $\beta=0.40$, respectively. The approximation of beta is done by taking the last $y$ value of each integral curve and divide it by the $x$ value ranging in $[\num{1.0e-5},\num{1.0e-4}]$, the adjusted beta is calculated using Equation~\eqref{2D Beta Adjust equation}. The points-axis corresponds to the 250 $x$ values we considered starting from $\num{1.0e-4}$ and ending in $\num{1.0e-5}$; that is, point $0$ corresponds to the value $x=\num{1.0e-4}$ and point $250$ corresponds to the value $x=\num{1.0e-5}$. The constants $c(y_0)=\ln(0.8)$ and $c(y_0)=1.2\ln(1.1)$ for the adjusted approximation correspond to Figure~\ref{fig:2D beta approx} $a)$ and $b)$, respectively. The approximations show an error that tends to decrease as we get closer to $x=0$ as depicted in the graphs. The value $\beta=0.40$ and $\beta=0.266$ are obtained by taking $a_0=15.0$, $a_2=5.0$, $b_0=1.0$ and $b_2=3.0$ and $a_0=15.0$, $a_2=6.0$, $b_0=1.0$ and $b_2=2.0$, respectively, in the vector field from Equation~\eqref{Normal Flow 2D}.

\begin{center}\vspace{1cm}
\includegraphics[width=1.1\linewidth]{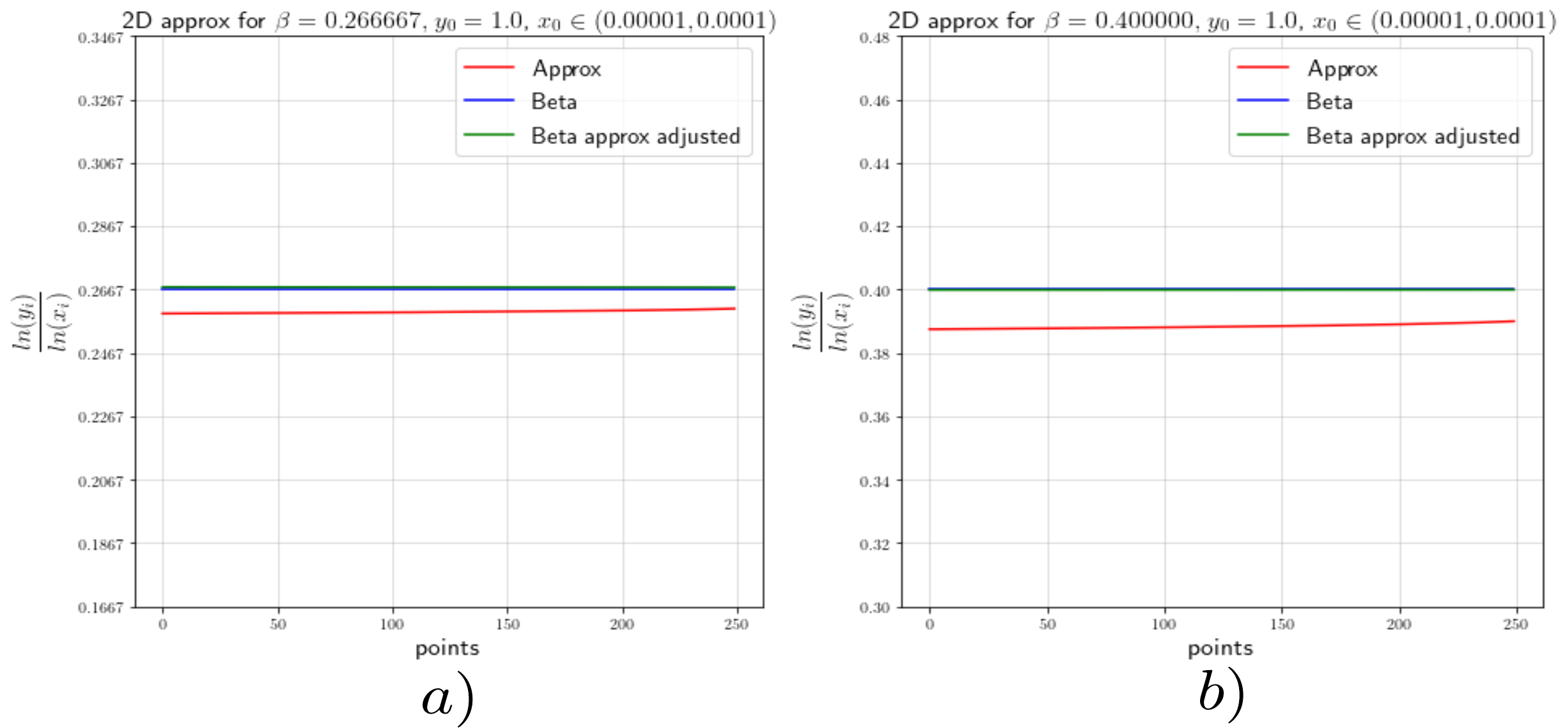}
\captionof{figure}{$2$-dimensional $\beta$ approximation.}
\label{fig:2D beta approx}
\end{center}\vspace{1cm}

We proceed in the same way to perform the numerical analysis for the exponent $\beta_2$. From the estimates obtained in \cite{Regular}, we can see that,

\begin{equation}
\beta_2\approx-\frac{\ln(x)}{\ln(t)},
\end{equation}

where $t$ is the flow time that it takes a point from the unstable leaf to hit the stable leaf and $x\in[\num{1.0e-5},\num{1.0e-4}]$. We will actually show that 

\begin{equation}\label{Beta2 Adjust equation 2D}
\beta_2=\frac{\ln(c(y_0))}{\ln(t)}-\frac{\ln(x)}{\ln(t)}+\mathcal{O}\Big(\frac{1}{\ln(t)}\Big).
\end{equation}


%

Figure~\ref{fig:2D beta2 approx} shows the approximation of the exponent of the tail of the return map with values $\beta_2=0.1333$ and $\beta_2=2.0$ which correspond to the values $\beta=0.400$ and $\beta=0.266$, respectively. For this case the constants for the adjusted approximations are $c(y_0)=\ln(0.06)$ and $c(y_0)=\ln(0.045)$ for Figure~\ref{fig:2D beta2 approx} $a)$ and $b)$, respectively. From the approximations we see, as before, that the error decreases as we approach to $x=0$.

\begin{center}\vspace{1cm}
\includegraphics[width=1.1\linewidth]{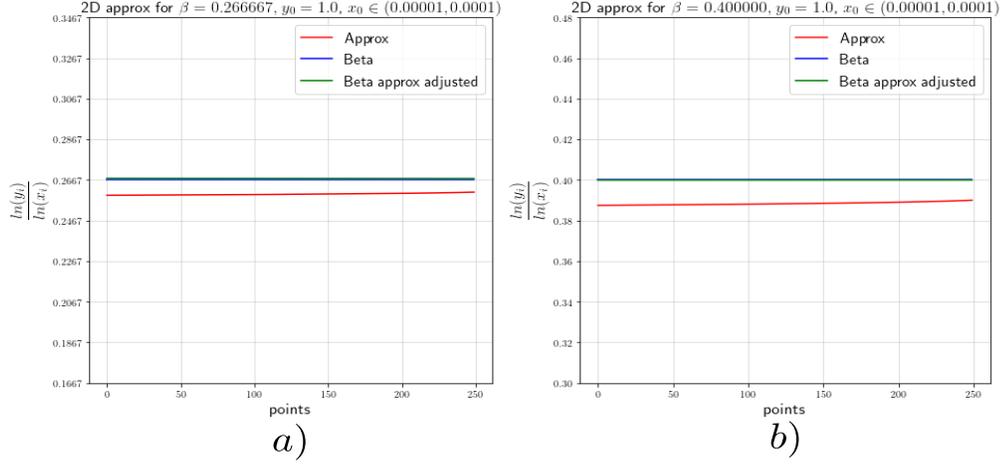}
\captionof{figure}{$2$-dimensional $\beta_2$ approximation.}
\label{fig:2D beta2 approx}
\end{center}\vspace{1cm}

\subsection{Numerics of the $3$-dimensional case}

In this subsection we will perform the numeric experiments for the $3$-dimensional models.
We will start with Neutral model 1. Recall that the Neutral model 1 was given by Equation~\eqref{Neutral model1}
\begin{equation}
\label{Neutral model1} 
\left( \begin{array}{ccc}
\dot{x}\\
\dot{y}\\
\dot{z}\\
\end{array}\right)=
N\left( \begin{array}{ccc}
x\\
y\\
z\\
\end{array}\right)=
\left( \begin{array}{ccc}
x(a_0x^2+a_1y^2+a_2z^2)\\
-\ell y\\
-z(b_0x^2+b_1y^2+b_2z^2)\\
\end{array}\right)+\mathcal{O}(4),  
\end{equation}
where $a_0,\;a_1,\;a_2,\;b_0,\;b_1,\;b_2$ and $\ell>0$ and $\Delta:=a_2b_0-a_0b_2\neq0$.Precise asymptotics are not available, but since $y(t)$ decreases exponentially fast, the same asymptotics as in \eqref{Normal Flow 2D} are expected, and our numerics indeed confirm this. Note that the strong stable direction of \eqref{Neutral model1} is still purely $y$-directed.


For the analysis of the Dulac map close to the neutral equilibrium of Equation~\eqref{Neutral model1}, we will perform the numerical analysis on $\fun{N_1}{\Sigma}{S}$. For the purpose of this work, we want to show with the numeric experiment that the $x$ and $z$ components behave like the $2$-dimensional model from the previous subsection regardless of the $y$ value. To perform the numerical analysis we will take different unstable leaves $W^u(x,y_0,z_0)$ and a stable leaf $W^s(1,0,0)$, where $y_0=1.0$ and $z_0=1.0$. Hence, like in the $2$-dimensional analysis we expect the numerical experiments to show us that,

\begin{equation}
\beta\approx\frac{\ln(z)}{\ln(x)},
\end{equation}

where $z$ is the last value of the integral curve with initial condition $(x,y_0,z_0)$ for $x\in[\num{1.0e-5},\num{1.0e-4}]$. Again, we will actually show that 

\begin{equation}\label{Beta Adjust equation 3D}
\beta=\frac{\ln(z)}{\ln(x)}-\frac{\ln(c(z_0))}{\ln(x)}+\mathcal{O}\Big(\frac{1}{\ln(x)}\Big).
\end{equation}


%

As before, we will use the Radau quadrature method and take 250 points for the values of $x$ starting from $\num{1.0e-4}$ and ending with $\num{1.0e-5}$; that is, point $0$ and point $250$ correspond to $x=\num{1.0e-4}$ and $x=\num{1.0e-5}$, respectively. We will consider the same values of $\beta$ we considered in the previous subsection. The approximation of $\beta$, corresponding to the red line in all figures, is done by taking the last $z$ value of each integral curve with initial condition $(x,y_0,z_0)$ and divide it by the $x$ value ranging in $[\num{1.0e-5},\num{1.0e-4}]$, the adjusted $\beta$, plotted in green in all figures, is calculated by using Equation~\eqref{Beta Adjust equation 3D}, and the theoretical value of $\beta$, corresponding to the blue graph in all figures, is obtained from the parameters $a_0,\; a_2,\;b_0$ and $b_2$ as before. The constants $c(z_0)$ were calculated using the least squares method. The constants  $c(z_0)=\ln(1.1)$ and $c(z_0)=\ln(1.06)$ correspond to Figure~\ref{fig:Neu1 beta approx} $a)$ and $b)$, respectively.

\begin{center}\vspace{1cm}
\includegraphics[width=1.1\linewidth]{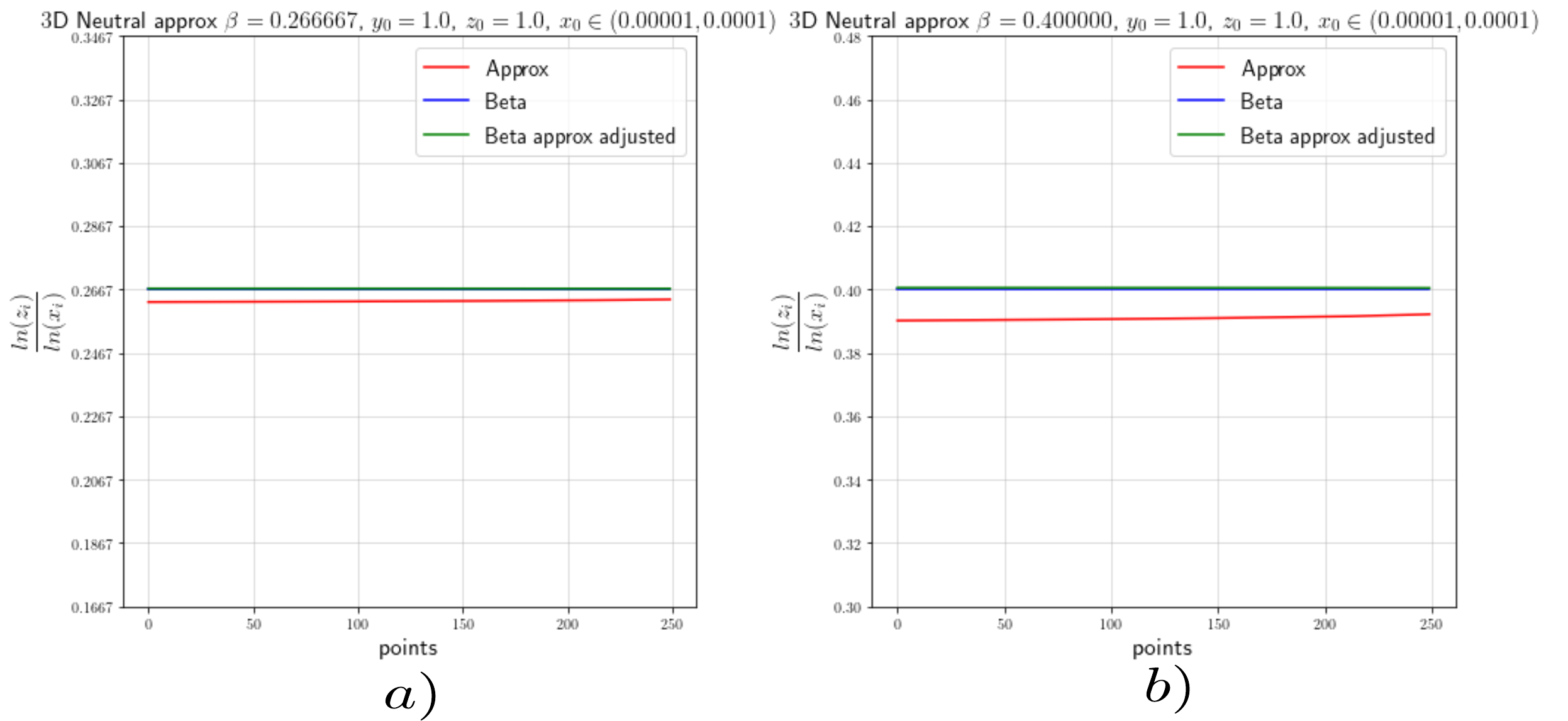}
\captionof{figure}{Neutral model 1 $\beta$ approximation.}
\label{fig:Neu1 beta approx}
\end{center}\vspace{1cm}

Next, we consider the Neutral model 2 given by Equation \eqref{Neutral model2}, with the parameters satisfying the usual constraints, and present the numerical results obtained by performing the same experiments we did for the Neutral model 1. We consider this form since it is no longer a skew product like the previous model and poses a new challenge to deduce its asymptotics and u decay of correlations.

\begin{equation}
\label{Neutral model2} 
\left( \begin{array}{ccc}
\dot{x}\\
\dot{y}\\
\dot{z}\\
\end{array}\right)=
G\left( \begin{array}{ccc}
x\\
y\\
z\\
\end{array}\right)=
\left( \begin{array}{ccc}
x(a_0x^2+a_2z^2)\\
-\ell y(1+c_0x^2+c_2z^2)\\
-z(b_0x^2+b_2z^2)\\
\end{array}\right)+\mathcal{O}(4).  
\end{equation}


%

Figure~\ref{fig:Neu2 beta approx} $a)$ and $b)$ show us the the numerical approximations of the Neutral model 2 for $\beta=0.40$ and $\beta=0.266$ with constants $c(z_0)=\ln(1.2)$ and $c(z_0)=\ln(1.08)$, respectively. We observe that the constants $c(z_0)$ and $c(y_0)$ from Equation~\eqref{Beta Adjust equation 3D} and Equation~\eqref{2D Beta Adjust equation}, respectively, are almost equal; that is, the $x$ and $z$ components of the $3$-dimensional model behaves like the $x$ and $y$ component of the $2$-dimensional model.

\begin{center}\vspace{1cm}
\includegraphics[width=1.1\linewidth]{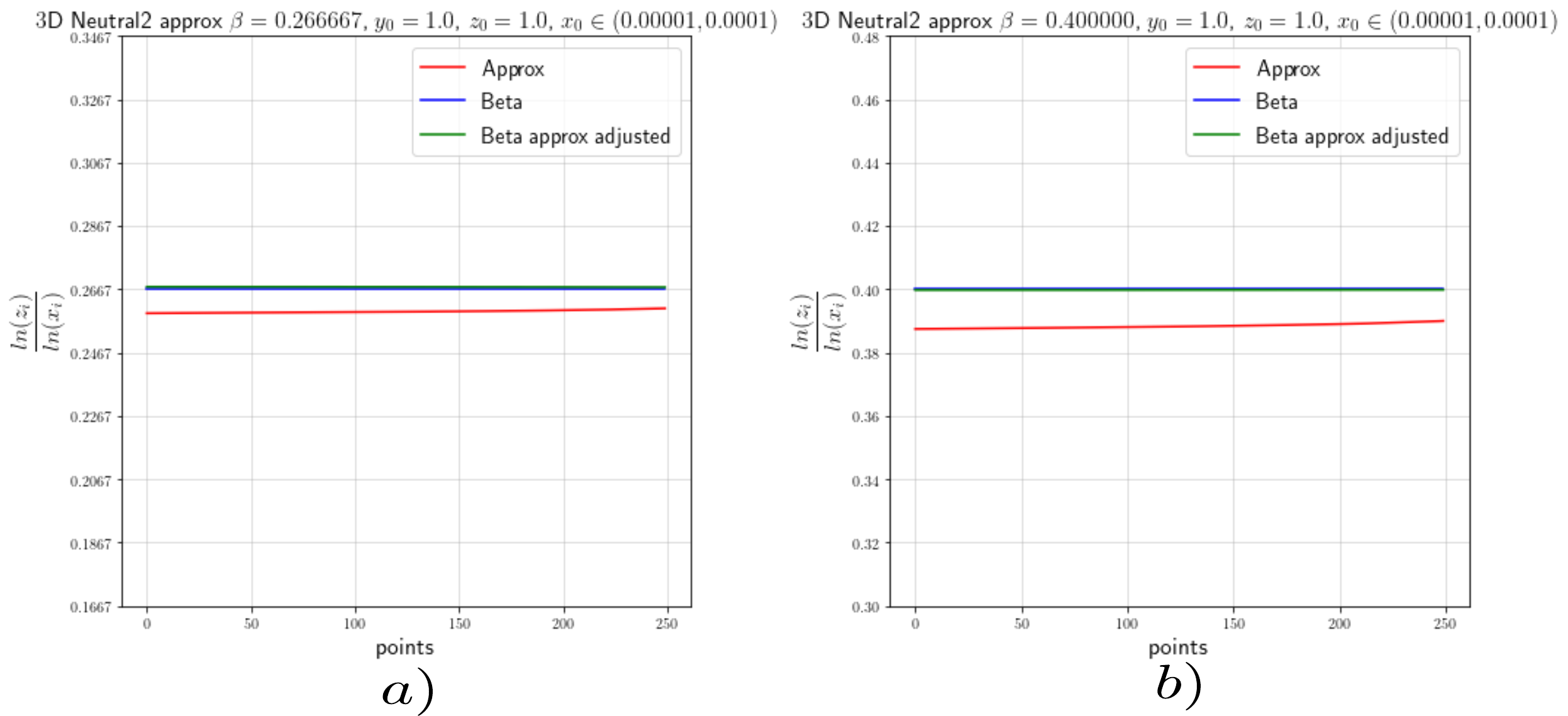}
\captionof{figure}{Neutral model 2 $\beta$ approximation.}
\label{fig:Neu2 beta approx}
\end{center}\vspace{1cm}

Until now we have performed the numerical experiments for the Neutral model 1 and 2 corresponding to the normal forms given by Equations~\eqref{Neutral model1} and \eqref{Neutral model2}, respectively. We saw there that for both models, the asymptotic behaviour of the $x$ and $z$ components are the same as the $x$ and $y$ component of the $2$-dimensional model. Recall that for this two models we had explicit formulas for the map $\fun{N}{\Sigma}{S}$ and hence for the modified Poincar\'e map $\fun{\Pneu}{\Sigma}{\Sigma}$. Our goal is to see whether the asymptotic behaviour of the Neutral model 3 given by Equation~\eqref{Neutral system} is similar to the other two models. Therefore, the numerical results obtained for the Neutral models 1 and 2 will be our reference and we will compare them to the numerical results obtained for the third model.


%

Figure~\ref{fig:Neu3 beta approx} $a)$ and $b)$ show us the the numerical approximations of the Neutral model 2 for $\beta=0.40$ and $\beta=0.266$ with constants $c(z_0)=\ln(1.12)$ and $c(z_0)=\ln(1.07)$, respectively. We observe again that the $x$ and $z$ components of the $3$-dimensional model behaves like the $x$ and $y$ component of the $2$-dimensional model. From this numerical experiments we can conclude that the behaviour of the map $N_1$ obtained by considering the neutral model 3 is asymptotically similar to the other two models; that is, the asymptotics of the Dulac map of the three models are similar.

\begin{center}\vspace{1cm}
\includegraphics[width=1.1\linewidth]{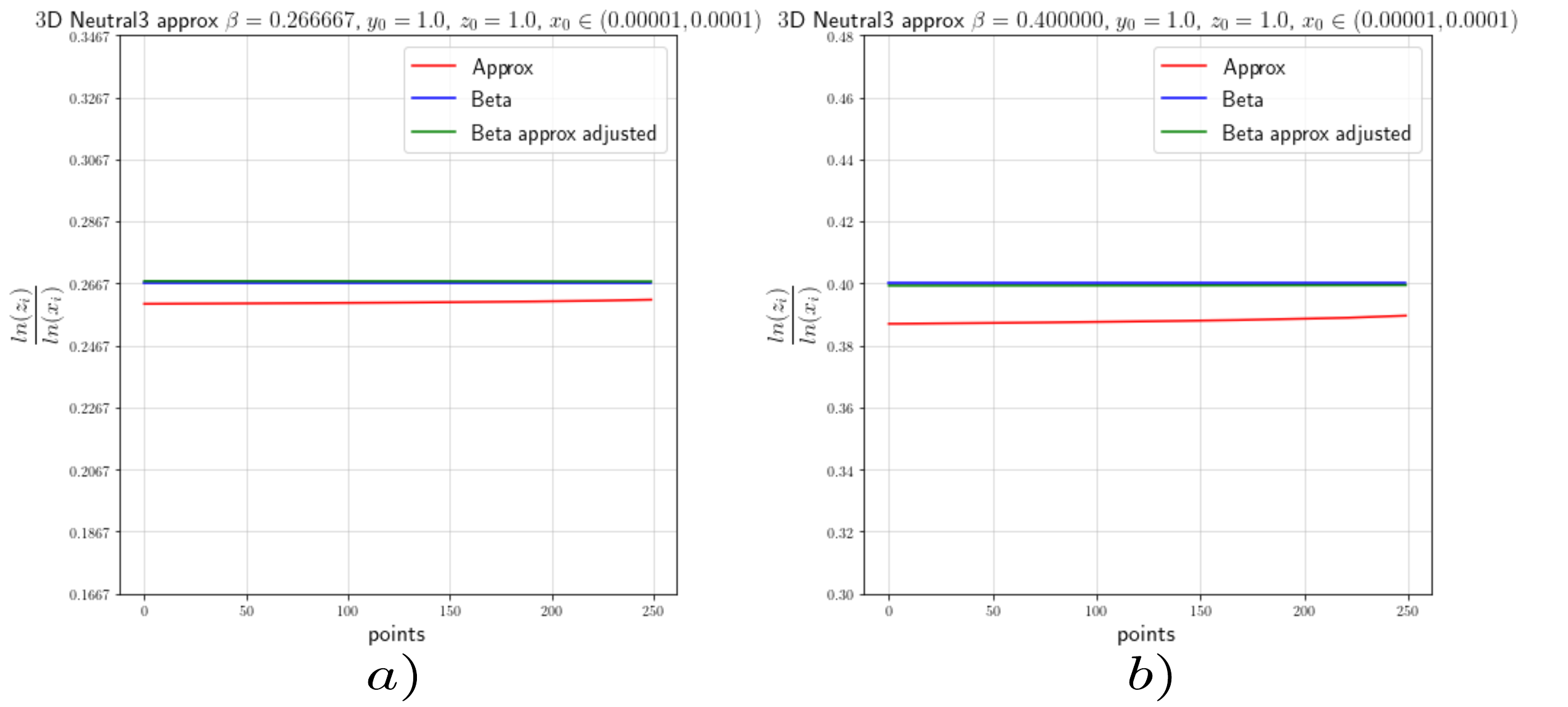}
\captionof{figure}{Neutral model 2 $\beta$ approximation.}
\label{fig:Neu3 beta approx}
\end{center}\vspace{1cm}

\subsection{Numerics of the tails of the return map in the \\$3$-dimensional case}

Next we will perform the numeric experiments for the $3$-dimensional models and see the approximations for the exponent of the decay of correlations; that is, for the exponent $\beta_2$. The general Neutral model or neutral model 3 is given by Equation~\eqref{Neutral system}, where $a_0,\;a_1,\;a_2,\;b_0,\;b_1,\;b_2,\;c_0,\;c_2$ and $\ell>0$ and $\Delta:=a_2b_0-a_0b_2\neq0$. Note that the neutral model 1 and the neutral model 2 are obtained from the general neutral model if we let $c_0,\;c_2=0$ and if we let $a_1,\;b_1=0$, respectively.

In the previous subsection we saw the numerical analysis on $\fun{N_1}{\Sigma}{S}$ and showed, with the numeric experimentation, that the $x$ and $z$ components behaves like the $2$-dimensional model. For the next numerical analysis, we will take an unstable leaf $W^u(x,y_0,z_0)$ and a stable leaf $W^s(1,0,0)$, where $y_0=1.0$ and $z_0=1.0$. From the estimates obtained in \cite{Regular} we can see that

\begin{equation}
\beta_2\approx-\frac{\ln(x)}{\ln(t)},
\end{equation}

where $t$ is the flow time that it takes a point from the unstable leaf to hit the stable leaf and $x\in[\num{1.0e-5},\num{1.0e-4}]$. We will actually show that 

\begin{equation}\label{Beta2 Adjust equation 3D}
\beta_2=\frac{\ln(c(z_0))}{\ln(t)}-\frac{\ln(x)}{\ln(t)}+\mathcal{O}\Big(\frac{1}{\ln(t)}\Big).
\end{equation}

As before, we will use the Radau quadrature method and take 50 points for the values of $x$ starting from $\num{1.0e-4}$ and ending with $\num{1.0e-5}$; that is, point $0$ and point $50$ correspond to $x=\num{1.0e-4}$ and $x=\num{1.0e-5}$, respectively. We will consider the same values of $\beta$ we considered in the previous subsection. The approximation of $\beta_2$, corresponding to the red line in all figures, is done by taking the $x$ value, ranging in $[\num{1.0e-5},\num{1.0e-4}]$, of each integral curve with initial condition $(x,y_0,z_0)$ and divide it by the flow time $t$, the adjusted approximation of $\beta_2$, shown in green in all figures, is calculated by using Equation~\eqref{Beta2 Adjust equation 3D}, and the theoretical value of $\beta_2$, corresponding to the blue graph in all figures, is obtained from the parameters $a_2$ and $b_2$; that is, $\beta_2=\frac{a_2+b_2}{2b_2}$.

We start considering the neutral model 1. Figure~\ref{fig:Neu1 beta2 approx} $a)$ shows the approximation for $\beta_2=1.333$ which corresponds to the case $\beta=0.40$, for the adjusted approximation the constant is $\ln(c(z_0))=0.06$ and $b)$ displays the approximation for $\beta_2=2.0$ corresponding to the case $\beta=0.266$ with adjustment constant $\ln(c(z_0))=0.05$.

\begin{center}\vspace{1cm}
\includegraphics[width=1.1\linewidth]{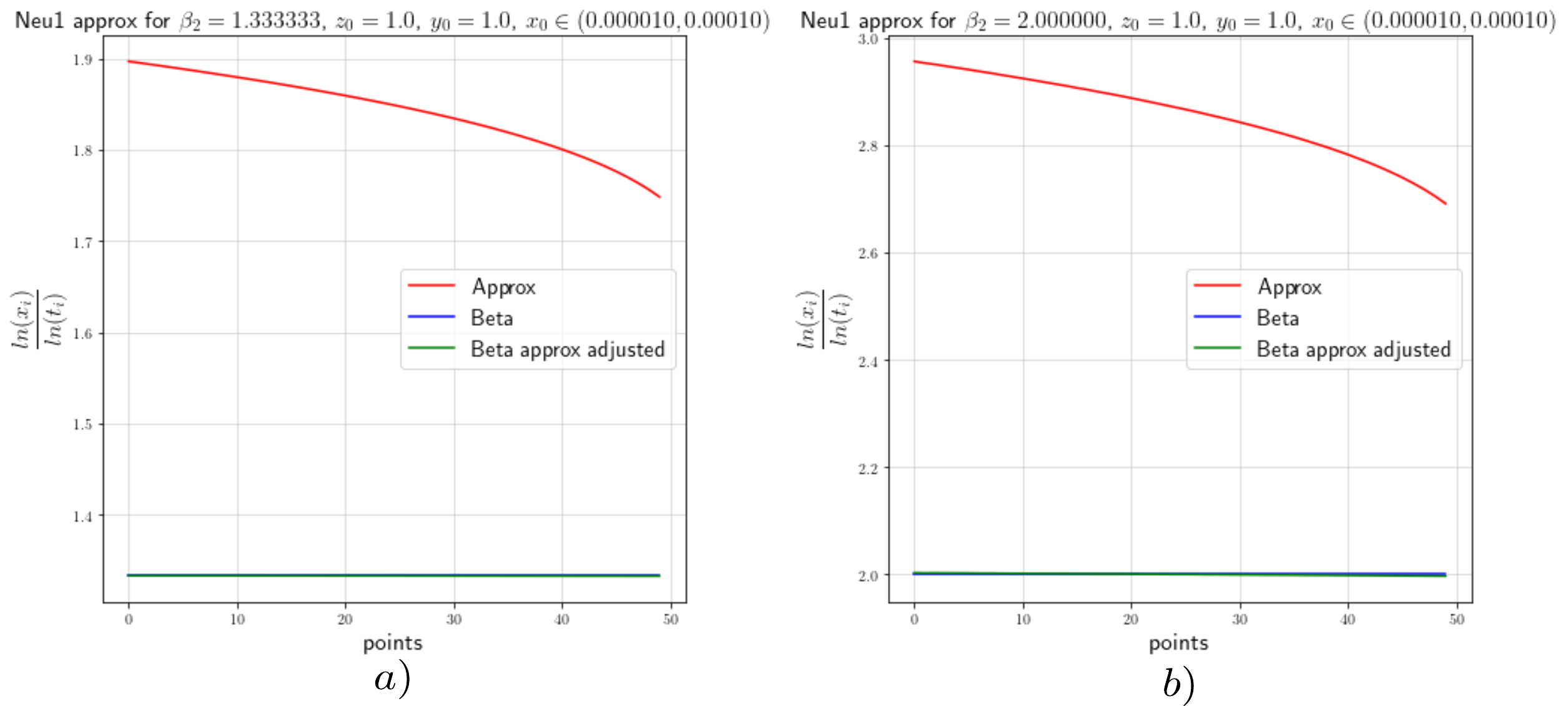}
\captionof{figure}{Neutral model 2 $\beta$ approximation.}
\label{fig:Neu1 beta2 approx}
\end{center}\vspace{1cm}


%
%
%

Next, we consider the neutral model 2. Figure~\ref{fig:Neu2 beta2 approx} $a)$ and $b)$ show the approximation for $\beta_2=1.333$ and $\beta_2=2.0$, respectively. Their adjusted approximation the constant are $\ln(c(z_0))=0.06$ and $\ln(c(z_0))=0.04$, respectively.

\begin{center}\vspace{1cm}
\includegraphics[width=1.1\linewidth]{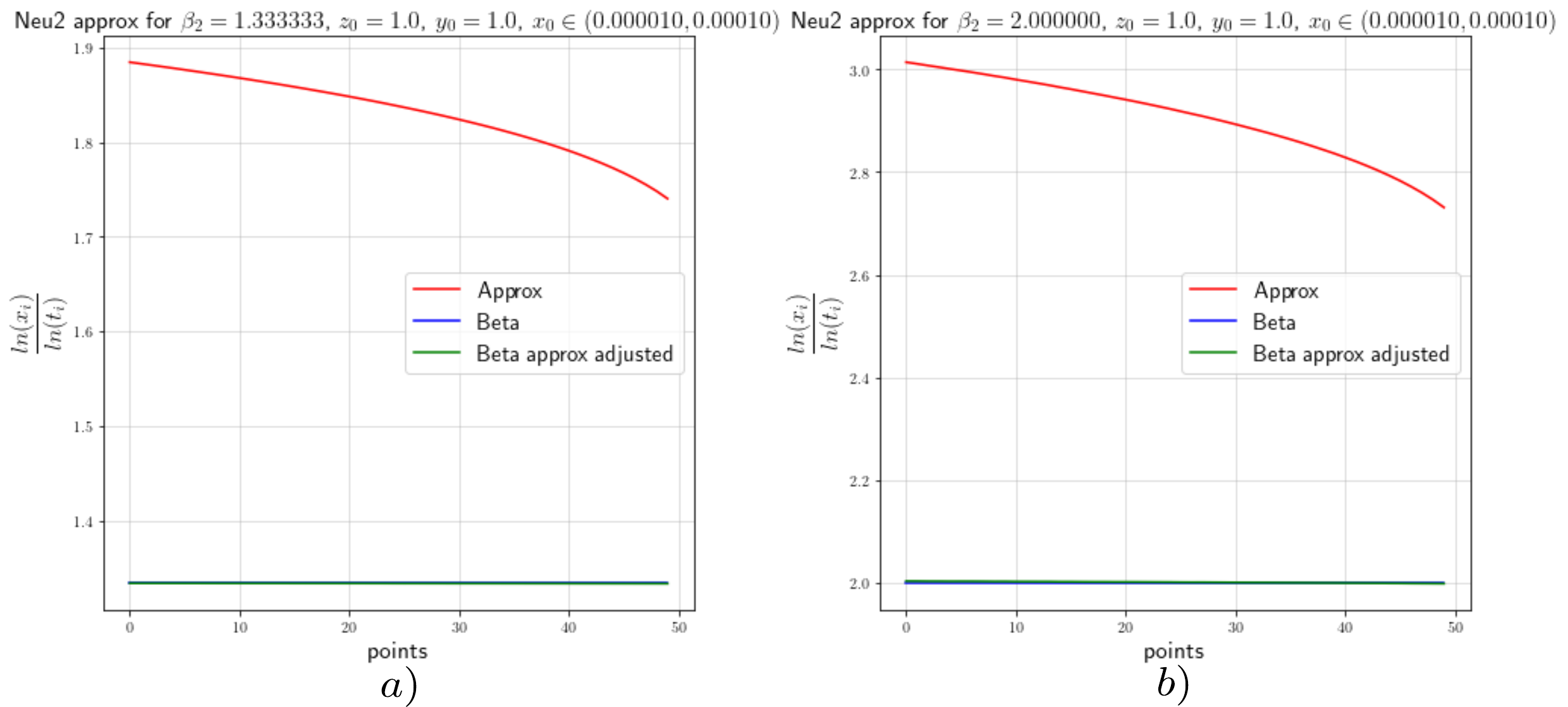}
\captionof{figure}{Neutral model 2 $\beta$ approximation.}
\label{fig:Neu2 beta2 approx}
\end{center}\vspace{1cm}

Finally, we consider the neutral model 3. Figure~\ref{fig:Neu3 beta2 approx} $a)$ and $b)$ show the approximation for $\beta_2=1.333$ and $\beta_2=2.0$, respectively. Their adjusted approximation the constant are $\ln(c(z_0))=0.06$ and $\ln(c(z_0))=0.04$, respectively.

\begin{center}\vspace{1cm}
\includegraphics[width=1.1\linewidth]{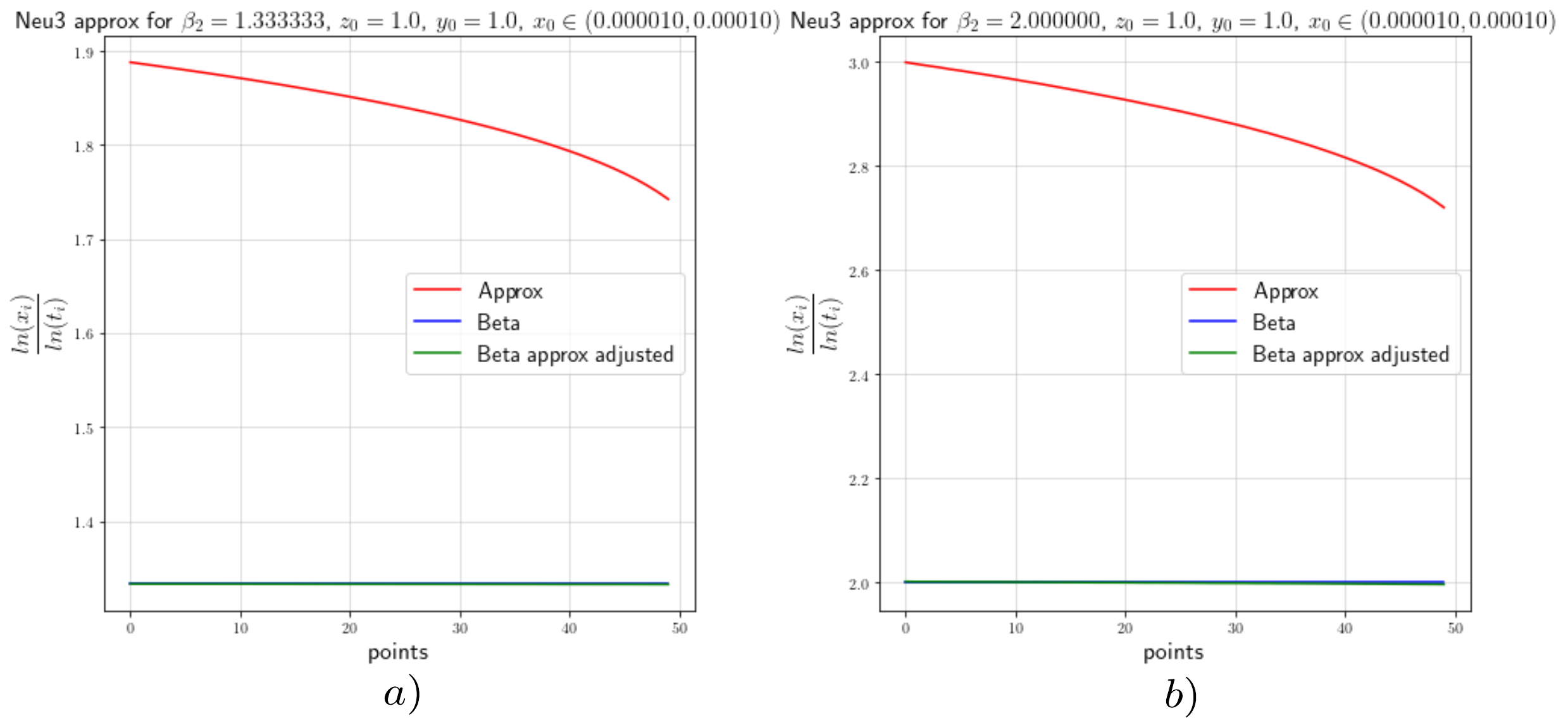}
\captionof{figure}{Neutral model 2 $\beta$ approximation.}
\label{fig:Neu3 beta2 approx}
\end{center}\vspace{1cm}

This numerical experiment has shown us a good approximation of the exponent of the decay of correlations for the 3 neutral models. From this we can deduce the same results concerning the decay of correlations, and obtaining Theorem~\ref{Poly decay Neutral geom Lorenz flow Thm} in its full generality.

\section{Programming specifications and data reproducibility}

This section introduces the technical specifications of the employed machine for conducting numerical analysis. It outlines the specific compiler, programming language, and numerical packages utilized in the study. Additionally, the subsequent subsections present the code implemented for both 2-dimensional and 3-dimensional analyses.

The numerical analysis was carried out using a Dell G3 series model 3579 laptop featuring a 15-inch display. The laptop is equipped with an 8th generation Intel Core i5 processor, coupled with a Mobile Intel HM370 chipset. The available RAM capacity is 8.0 GB. The machine boasts two graphic cards: an integrated Intel UHD Graphics 630 and a discrete NVIDIA GeForce GTX 1050. The operating system employed is Windows 11 Home (64-bit).

For the analysis, the compiler utilized was Visual Studio Code version 1.81.1, while Python version 3.9.12 served as the programming language. The implemented code made use of the following packages: scipy 1.9.3, numpy 1.21.5, and matplotlib 3.5.1.

\subsection{Code for the 2-dimensional case}
\begin{verbatim}
"""This script provides the ploting of a the two
dimensionl system of ode used in my PhD studies

Author: Hector Homero Canales Farias
"""
# %%
import numpy as np
from scipy.integrate import solve_ivp
import matplotlib.pyplot as plt

# Parameters neutral system
a0, a2 = 15.0, 6.0
b0, b2 = 1.0, 2.0
beta0 = (a0+b0) / (2*a0)
beta2 = (a2+b2) / (2*b2)
beta = beta0 / beta2

# plt.rcParams["figure.figsize"] = (7,7)
plt.rcParams['text.usetex'] = True


# auxiliary functio for plotting
def beta_plot(x, value):
    return np.full(len(x), value)


# function to find constants for the adjustment of approx.
def find_constanst(num, dem, goal):
    B = np.log(num)/np.log(dem) - goal
    A = np.matrix(1 / np.log(dem)).T
    return np.linalg.lstsq(A, B, rcond=None)[0][0]


# Normal equation for solve_ivp with condition stopping at x=1
def Neutral_stop(time, y, a_0, a_2, b_0, b_2):

    if abs(y[0]) <= 1:

        return np.array([
            y[0]*(a_0*(y[0]**2) + a_2*(y[1]**2)),
            -y[1]*(b_0*(y[0]**2) + b_2*(y[1]**2)),
        ])

    else:

        return np.array([0.0, 0.0])


# here the ode will be solve and list for the last y-values are recorded
t1 = [0.0, 10500.0]
dt1 = 0.01
teval1 = np.arange(0.0, 10500.0, dt1)
a = np.linspace(0.0001, 0.00001, 250)
y0 = 1.0
x_initial1 = [np.array([a[i], y0]) for i in range(a.shape[0])]

x1 = []
y1 = []
states_orbits1 = []
temp1 = []
x_last1 = []
y_last1 = []  # y-value when x hist 1

for i in range(len(x_initial1)):

    u1 = solve_ivp(
         Neutral_stop, t1, x_initial1[i],
         args=(a0, a2, b0, b2), method='Radau',
         t_eval=teval1
    )
    x_last1.append(u1.y[0, -1])
    y_last1.append(u1.y[1, -1])
    temp1.append(u1.t)
    x1.append(u1.y[0, :])
    y1.append(u1.y[1, :])
    states_orbits1.append(u1.y)

# Since different conditions require different times for x to hit 1
# record the first index when x hits 1 in a list
index = []
for i in range(len(x1)):
    indices = np.where(x1[i] >= 1.0)
    index.append(indices[0][0])

# list of the times for those exact indices
t_last = []
for i in range(len(temp1)):
    t_last.append(temp1[i][index[i]])

# list of the y-values at the given indices
y_last = []
for i in range(len(index)):
    y_last.append(y1[i][index[i]])

# define the adjustment constants
const0 = find_constanst(y_last, t_last, -beta0)
const2 = find_constanst(a, t_last, -beta2)
const = find_constanst(y_last, a, beta)

# define beta cosntants
beta_asymp = np.log(y_last)/np.log(a)
beta2_asymp = -np.log(a)/np.log(t_last)
beta0_asymp = -np.log(y_last)/np.log(t_last)

# define beta constants adjusted
beta_asymp_adjust = beta_asymp - const/np.log(a)
beta2_asymp_adjust = beta2_asymp + const2/np.log(t_last)
beta0_asymp_adjust = beta0_asymp + const0/np.log(t_last)

# plot of the beta approx
title = (
    r'2D approx for $\beta={:0.6f},$ $y_0={:0.1f}$, $ x_0\in({:0.5f},{:0.4f})$'
)
title_2 = (
    r'2D approx for $\beta_2={:0.6f},$ $y_0={:0.1f}$,'
    r'$ x_0\in({:0.5f},{:0.4f})$'
)
y_label1 = (
    r'$\displaystyle\frac{ln(y_i)}{ln(x_i)}$ '
)
y_label2 = (
    r'$\displaystyle\frac{ln(x_i)}{ln(t_i)}$ '
)
fig1 = plt.figure(figsize=(7, 7))
ax1 = plt.gca()

ax1.plot(range(len(beta_asymp)),
         beta_asymp, color="red",
         label='Approx')
ax1.plot(range(len(beta_asymp)),
         beta_plot(beta_asymp, beta),
         color='blue',
         label='Beta')
ax1.plot(range(len(beta_asymp_adjust)),
         beta_asymp_adjust, color='green',
         label='Beta approx adjusted')
ax1.set_yticks(np.arange(beta-0.1, max(beta_asymp)+0.1, 0.02))
ax1.set_xlabel("points", fontsize=16)
ax1.set_ylabel(y_label1, fontsize=16)
ax1.set_title(title.format(
    beta, y0, a[-1], a[0]), fontsize=16)
ax1.grid(alpha=0.5)
ax1.legend(prop={'size': 15})
plt.show()

# plot of the beta2 approx
fig2 = plt.figure(figsize=(7, 7))
ax2 = plt.gca()

ax2.plot(range(len(beta2_asymp)),
         beta2_asymp, color="red",
         label='Approx')
ax2.plot(range(len(beta2_asymp)),
         beta_plot(beta2_asymp, beta2),
         color='blue', label='Beta')
ax2.plot(range(len(beta2_asymp_adjust)),
         beta2_asymp_adjust, color='green',
         label='Beta approx adjusted')
ax2.set_xlabel("points", fontsize=16)
ax2.set_ylabel(y_label2, fontsize=16)
ax2.set_title(title.format(
    beta2, y0, a[-1], a[0]), fontsize=16)
ax2.grid(alpha=0.5)
ax2.legend(prop={'size': 15})
plt.show()
\end{verbatim}

\subsection{Code for the 3-dimensional case}

\begin{verbatim}
"""This script provides the ploting of a the three
dimensionl system of ode used in my PhD studies

Author: Hector Homero Canales Farias
"""
# %%
import numpy as np
from scipy.integrate import solve_ivp
import matplotlib.pyplot as plt

# Parameters neutral system
a0, a1, a2 = 15.0, 4.0, 6.0
b0, b1, b2 = 1.0, 2.0, 2.0
c0, c2 = 1.0, 3.0
l = 8
c1 = 2.8
beta0 = (a0+b0) / (2*a0)
beta2 = (a2+b2) / (2*b2)
beta = beta0 / beta2

# plt.rcParams["figure.figsize"] = (7,7)
plt.rcParams['text.usetex'] = True


# auxiliary functio for plotting
def beta_plot(x, value):
    return np.full(len(x), value)


# function to find constants for the adjustment of approx.
def find_constanst(num, dem, goal):
    B = np.log(num)/np.log(dem) - goal
    A = np.matrix(1 / np.log(dem)).T
    return np.linalg.lstsq(A, B, rcond=None)[0][0]


# Normal equation for solve_ivp with stopping
def Neutral_stop(time, y, a_0, a_1, a_2, b_0, b_1, b_2, l_cst):

    if abs(y[0]) <= 1:

        return np.array([
            y[0]*(a_0*(y[0]**2) + a_1*(y[1]**2) + a_2*(y[2]**2)),
            -l_cst*y[1],
            -y[2]*(b_0*(y[0]**2) + b_1*(y[1]**2) + b_2*(y[2]**2)),
        ])

    else:

        return np.array([0.0, 0.0, 0.0])


# Normal2 for solve_ivp with stopping
def Neutral2_stop(time, y, a_0, a_2, b_0, b_2, c_0, c_2, l_cst):

    if abs(y[0]) <= 1:

        return np.array([
            y[0]*(a_0*(y[0]**2) + a_2*(y[2]**2)),
            -l_cst*y[1]*(1+c_0*(y[0]**2)+c_2*(y[2]**2)),
            -y[2]*(b_0*(y[0]**2) + b_2*(y[2]**2)),
        ])

    else:

        return np.array([0.0, 0.0, 0.0])


# Normal3 for solve_ivp with stopping
def Neutral3_stop(time, y, a_0, a_1, a_2, b_0, b_1, b_2, c_0, c_2, l_cst):

    if abs(y[0]) <= 1:

        return np.array([
            y[0]*(a_0*(y[0]**2) + a_1*(y[1]**2) + a_2*(y[2]**2)),
            -l_cst*y[1]*(1+c_0*(y[0]**2)+c_2*(y[2]**2)),
            -y[2]*(b_0*(y[0]**2) + b_1*(y[1]**2) + b_2*(y[2]**2)),
        ])

    else:
        return np.array([0.0, 0.0, 0.0])


# Neutral Case model 1
t1 = [0.0, 3500.0]
dt1 = 0.01
teval1 = np.arange(0.0, 3500.0, dt1)
y_init1 = 1.0
z_init1 = 1.0
a_init1 = np.linspace(0.0001, 0.00001, 50)
x_initial1 = [
    np.array([a_init1[i], y_init1, z_init1]) for i in range(a_init1.shape[0])
]

x1 = []
y1 = []
z1 = []
states_orbits1 = []
temp1 = []

for i in range(len(x_initial1)):
    u1 = solve_ivp(
        Neutral_stop, t1, x_initial1[i],
        args=(a0, a1, a2, b0, b1, b2, l), method='Radau',
        t_eval=teval1)
    temp1.append(u1.t)
    x1.append(u1.y[0, :])
    y1.append(u1.y[1, :])
    z1.append(u1.y[2, :])
    states_orbits1.append(u1.y)

# record the first index when x hits 1 in a list
# record corresponding times and z-values
index1 = []
for i in range(len(x1)):
    indices = np.where(x1[i] >= 1.0)
    index1.append(indices[0][0])

t_last1 = []
for i in range(len(temp1)):
    t_last1.append(temp1[i][index1[i]])

z_last1 = []
for i in range(len(index1)):
    z_last1.append(z1[i][index1[i]])

# define adjustment constants
cst2_Neu1 = find_constanst(a_init1, t_last1, -beta2)
cst_Neu1 = find_constanst(z_last1, a_init1, beta)

# define beta constants and cst with adjustment
beta_asympNeu1 = np.log(z_last1)/np.log(a_init1)
beta2_asympNeu1 = -np.log(a_init1)/np.log(t_last1)
beta_asymp_adjustNeu1 = beta_asympNeu1 - cst_Neu1/np.log(a_init1)
beta2_asymp_adjustNeu1 = beta2_asympNeu1 + cst2_Neu1/np.log(t_last1)

# plot beta constant
title_beta_neu1 = (
    r'Neu1 $\beta={:0.6f}$, $y_0={:0.1f}$, $z_0={:0.1f}$,'
    r'$x_0\in({:0.6f},{:0.5f})$'
)
ylabel_beta_neu = (
    r'$\displaystyle\frac{ln(z_i)}{ln(x_i)}$ '
)
fig1_1 = plt.figure(figsize=(7, 7))
ax1_1 = fig1_1.add_subplot()
ax1_1.plot(range(len(beta_asympNeu1)), beta_asympNeu1,
           color="red", label='Approx')
ax1_1.plot(range(len(beta_asympNeu1)), beta_plot(beta_asympNeu1, beta),
           color='blue', label='Beta')
ax1_1.plot(range(len(beta_asymp_adjustNeu1)), beta_asymp_adjustNeu1,
           color='green', label='Beta approx adjusted')
ax1_1.set_xlabel("points", fontsize=16)
ax1_1.set_ylabel(ylabel_beta_neu, fontsize=16)
ax1_1.set_title(
    title_beta_neu1.format(
        beta, y_init1,
        z_init1, a_init1[-1],
        a_init1[0]
    ),
    fontsize=16
)
ax1_1.grid(alpha=0.5)
ax1_1.legend(prop={'size': 15})
plt.show()

# plot beta2 approx
title_beta2_neu1 = (
    r'Neu1 approx for $\beta_2={:0.6f}$,  $z_0={:0.1f}$, $y_0={:0.1f}$, $'
    r'x_0\in({:0.6f},{:0.5f})$'
)
ylabel_beta2_neu = (
    r'$\displaystyle\frac{ln(x_i)}{ln(t_i)}$ '
)
fig1_2 = plt.figure(figsize=(7, 7))
ax1_2 = plt.gca()

ax1_2.plot(range(len(beta2_asympNeu1)), beta2_asympNeu1,
           color="red", label='Approx')
ax1_2.plot(range(len(beta2_asympNeu1)), beta_plot(beta2_asympNeu1, beta2),
           color='blue', label='Beta')
ax1_2.plot(range(len(beta2_asymp_adjustNeu1)), beta2_asymp_adjustNeu1,
           color='green', label='Beta approx adjusted')
ax1_2.set_xlabel("points", fontsize=16)
ax1_2.set_ylabel(ylabel_beta2_neu, fontsize=16)
ax1_2.set_title(
    title_beta2_neu1.format(
        beta2, y_init1, z_init1,
        a_init1[-1], a_init1[0]
    ),
    fontsize=16)
ax1_2.grid(alpha=0.5)
ax1_2.legend(prop={'size': 15})
plt.show()

# Neutral model 2
t2 = [0.0, 3500.0]
dt2 = 0.01
teval2 = np.arange(0.0, 3500.0, dt2)
y_init2 = 1.0
z_init2 = 1.0
a_init2 = np.linspace(0.0001, 0.00001, 50)
x_initial2 = [
    np.array([a_init2[i], y_init2, z_init2]) for i in range(a_init2.shape[0])
]

x2 = []
y2 = []
z2 = []
states_orbits2 = []
temp2 = []

for i in range(len(x_initial2)):

    u4 = solve_ivp(Neutral2_stop, t2, x_initial2[i],
                   args=(a0, a2, b0, b2, c0, c2, l),
                   method='Radau', t_eval=teval2)
    temp2.append(u4.t)
    x2.append(u4.y[0, :])
    y2.append(u4.y[1, :])
    z2.append(u4.y[2, :])
    states_orbits2.append(u4.y)

# record the first index when x hits 1 in a list
# record corresponding times and z-values
index2 = []
for i in range(len(x2)):
    indices = np.where(x2[i] >= 1.0)
    index2.append(indices[0][0])

t_last2 = []
for i in range(len(temp2)):
    t_last2.append(temp2[i][index2[i]])

z_last2 = []
for i in range(len(index2)):
    z_last2.append(z2[i][index2[i]])

# define adjustment constants
cst2_Neu2 = find_constanst(a_init2, t_last2, -beta2)
cst_Neu2 = find_constanst(z_last2, a_init2, beta)

# define beta constants and cst with adjustment
beta_asympNeu2 = np.log(z_last2)/np.log(a_init2)
beta2_asympNeu2 = -np.log(a_init2)/np.log(t_last2)
beta_asymp_adjustNeu2 = beta_asympNeu2 - cst_Neu2/np.log(a_init2)
beta2_asymp_adjustNeu2 = beta2_asympNeu2 + cst2_Neu2/np.log(t_last2)

# plot beta constant
title_beta_neu2 = (
    r'Neu2 approx $\beta={:0.6f}$, $y_0={:0.1f}$, $z_0={:0.1f}$,'
    r'$x_0\in({:0.6f},{:0.5f})$'
)
fig2_1 = plt.figure(figsize=(7, 7))
ax2_1 = fig2_1.add_subplot()
ax2_1.plot(range(len(beta_asympNeu2)), beta_asympNeu2,
           color="red", label='Approx')
ax2_1.plot(range(len(beta_asympNeu2)), beta_plot(beta_asympNeu2, beta),
           color='blue', label='Beta')
ax2_1.plot(range(len(beta_asymp_adjustNeu2)), beta_asymp_adjustNeu2,
           color='green', label='Beta approx adjusted')
ax2_1.set_xlabel("points", fontsize=16)
ax2_1.set_ylabel(ylabel_beta_neu, fontsize=16)
ax2_1.set_title(
    title_beta_neu2.format(
        beta, y_init2,
        z_init2, a_init2[-1],
        a_init2[0]
    ),
    fontsize=16
)
ax2_1.grid(alpha=0.5)
ax2_1.legend(prop={'size': 15})
plt.show()

# plot beta2 approx
title_beta2_neu2 = (
    r'Neu2 approx for $\beta_2={:0.6f}$,  $z_0={:0.1f}$, $y_0={:0.1f}$, $'
    r'x_0\in({:0.6f},{:0.5f})$'
)
fig2_2 = plt.figure(figsize=(7, 7))
ax2_2 = plt.gca()

ax2_2.plot(range(len(beta2_asympNeu2)), beta2_asympNeu2,
           color="red", label='Approx')
ax2_2.plot(range(len(beta2_asympNeu2)), beta_plot(beta2_asympNeu2, beta2),
           color='blue', label='Beta')
ax2_2.plot(range(len(beta2_asymp_adjustNeu2)), beta2_asymp_adjustNeu2,
           color='green', label='Beta approx adjusted')
ax2_2.set_xlabel("points", fontsize=16)
ax2_2.set_ylabel(ylabel_beta2_neu, fontsize=16)
ax2_2.set_title(
    title_beta2_neu2.format(
        beta2, y_init2, z_init2,
        a_init2[-1], a_init2[0]
    ),
    fontsize=16)
ax2_2.grid(alpha=0.5)
ax2_2.legend(prop={'size': 15})
plt.show()

# Neutral model 3
t3 = [0.0, 3500.0]
dt3 = 0.01
teval3 = np.arange(0.0, 3500.0, dt3)
y_init3 = 1.0
z_init3 = 1.0
a_init3 = np.linspace(0.0001, 0.00001, 50)
x_initial3 = [
    np.array([a_init3[i], y_init3, z_init3]) for i in range(a_init3.shape[0])
]

x3 = []
y3 = []
z3 = []
states_orbits3 = []
temp3 = []

for i in range(len(x_initial3)):

    u3 = solve_ivp(Neutral3_stop, t3, x_initial3[i],
                   args=(a0, a1, a2, b0, b1, b2, c0, c2, l),
                   method='Radau', t_eval=teval3)
    temp3.append(u3.t)
    x3.append(u3.y[0, :])
    y3.append(u3.y[1, :])
    z3.append(u3.y[2, :])
    states_orbits3.append(u3.y)

# record the first index when x hits 1 in a list
# record corresponding times and z-values
index3 = []
for i in range(len(x3)):
    indices = np.where(x3[i] >= 1.0)
    index3.append(indices[0][0])

t_last3 = []
for i in range(len(temp3)):
    t_last3.append(temp3[i][index3[i]])

z_last3 = []
for i in range(len(index3)):
    z_last3.append(z3[i][index3[i]])

# define adjustment constants
cst2_Neu3 = find_constanst(a_init3, t_last3, -beta2)
cst_Neu3 = find_constanst(z_last3, a_init3, beta)

# define beta constants and cst with adjustment
beta_asympNeu3 = np.log(z_last3)/np.log(a_init3)
beta2_asympNeu3 = -np.log(a_init3)/np.log(t_last3)
beta_asymp_adjustNeu3 = beta_asympNeu3 - cst_Neu3/np.log(a_init3)
beta2_asymp_adjustNeu3 = beta2_asympNeu3 + cst2_Neu3/np.log(t_last3)

# plot beta constant
title_beta_neu3 = (
    r'Neu3 approx $\beta={:0.6f}$, $y_0={:0.1f}$, $z_0={:0.1f}$,'
    r'$x_0\in({:0.6f},{:0.5f})$'
)
fig3_1 = plt.figure(figsize=(7, 7))
ax3_1 = fig2_1.add_subplot()
ax3_1.plot(range(len(beta_asympNeu3)), beta_asympNeu3,
           color="red", label='Approx')
ax3_1.plot(range(len(beta_asympNeu3)), beta_plot(beta_asympNeu3, beta),
           color='blue', label='Beta')
ax3_1.plot(range(len(beta_asymp_adjustNeu3)), beta_asymp_adjustNeu3,
           color='green', label='Beta approx adjusted')
ax3_1.set_xlabel("points", fontsize=16)
ax3_1.set_ylabel(ylabel_beta_neu, fontsize=16)
ax3_1.set_title(
    title_beta_neu3.format(
        beta, y_init3,
        z_init3, a_init3[-1],
        a_init3[0]
    ),
    fontsize=16
)
ax3_1.grid(alpha=0.5)
ax3_1.legend(prop={'size': 15})
plt.show()

# plot beta2 approx
title_beta2_neu3 = (
    r'Neu3 approx for $\beta_2={:0.6f}$,  $z_0={:0.1f}$, $y_0={:0.1f}$, $'
    r'x_0\in({:0.6f},{:0.5f})$'
)
fig3_2 = plt.figure(figsize=(7, 7))
ax3_2 = plt.gca()

ax3_2.plot(range(len(beta2_asympNeu3)), beta2_asympNeu3,
           color="red", label='Approx')
ax3_2.plot(range(len(beta2_asympNeu3)), beta_plot(beta2_asympNeu3, beta2),
           color='blue', label='Beta')
ax3_2.plot(range(len(beta2_asymp_adjustNeu3)), beta2_asymp_adjustNeu3,
           color='green', label='Beta approx adjusted')
ax3_2.set_xlabel("points", fontsize=16)
ax3_2.set_ylabel(ylabel_beta2_neu, fontsize=16)
ax3_2.set_title(
    title_beta2_neu3.format(
        beta2, y_init3, z_init3,
        a_init3[-1], a_init3[0]
    ),
    fontsize=16)
ax3_2.grid(alpha=0.5)
ax3_2.legend(prop={'size': 15})
plt.show()
\end{verbatim}

\end{document}